\documentclass[11pt,a4paper]{article}
\usepackage{amssymb,amsmath}

\newtheorem{stat}{Statement}[section]
\newtheorem{prop}[stat]{Proposition}
\newtheorem{cor}[stat]{Corollary}
\newtheorem{thm}[stat]{Theorem}
\newtheorem{lem}[stat]{Lemma}

\newtheorem{rem}[stat]{Remark}

\newtheorem{defi}[stat]{Definition}

\newtheorem{nota}[stat]{Notation}

\numberwithin{equation}{section}

\newenvironment{proof}{\noindent {\bf Proof}.}{\hfill$\triangle$}

\def\Section{\setcounter{equation}{0}\section}

\def\be{\begin{eqnarray}}
\def\ee{\end{eqnarray}}

\begin{document}
\title{Estimates for the density of a nonlinear Landau process}
\author{H\'el\`ene Gu\'erin$^1$, Sylvie M\'el\'eard$^2$ and Eulalia Nualart$^3$}

\footnotetext[1]{IRMAR, Universit\'e Rennes 1, Campus de Beaulieu,
35042 Rennes, France. E-mail: helene.guerin@univ-rennes1.fr}

\footnotetext[2]{MODAL'X, Universit\'e Paris 10, 200 av. de la
R\'epublique, 92000 Nanterre, France. E-mail:
sylvie.meleard@u-paris10.fr}

\footnotetext[3]{Institut Galil\'ee, Universit\'e Paris 13, av.
J.-B. Cl\'ement, 93430 Villetaneuse, France. E-mail:
nualart@math.univ-paris13.fr} \maketitle

\begin{abstract}
The aim of this paper is to obtain estimates for the density of
the law of a specific nonlinear diffusion process at any positive
bounded time. This process is issued from kinetic theory and is
called Landau process, by analogy with the associated
deterministic Fokker-Planck-Landau equation. It is not Markovian,
its coefficients are not bounded and the diffusion matrix is
degenerate. Nevertheless, the specific form of the diffusion
matrix and the nonlinearity imply the non-degeneracy of the
Malliavin matrix and then the existence and smoothness of the
density. In order to obtain a lower bound for the density, the
known results do not apply. However, our approach follows the
main idea consisting in discretizing the interval time and
developing a recursive method. To this aim, we prove and use
refined results on conditional Malliavin calculus. The lower bound
implies the positivity of the solution of the Landau equation, and
partially answers to an analytical conjecture. We also obtain an
upper bound for the density, which again leads to an unusual
estimate due to the bad behavior of the coefficients.
\end{abstract}
\vskip 1,5cm {\it \noindent AMS 2000 subject classifications:}
Primary: 60H30, 60H07; Secondary: 82C31,  82C40. \vskip 10pt
\noindent {\it Key words and phrases}. Conditional Malliavin
calculus, density estimates, nonlinear Landau process, unbounded
coefficients, Fokker-Planck-Landau equation. \vskip 5cm \pagebreak

\Section{Introduction} In this paper, we consider a nonlinear
diffusion process issued from kinetic theory and called Landau
process, by analogy with the associated deterministic Landau
equation. This process is defined as the solution of a nonlinear
stochastic differential equation driven by a space-time white
noise. Its coefficients are obtained from the Landau equation. In
particular, they are not bounded and the diffusion matrix is
degenerate. Nevertheless, Gu\'erin \cite{Guerin:02} uses the
nonlinearity of the equation and the specific form of the
diffusion matrix to prove the existence and smoothness of the
density of the law of this process at each finite time. This
implies in particular the existence of a smooth solution to the
nonlinear partial differential Landau equation.

The aim of this paper is
 to obtain lower and upper bounds for this density. The bad
 behavior of the coefficients of the stochastic differential
 equation makes the problem unusual. In particular,
 the methods introduced by Kusuoka and
Stroock~\cite{Kusuoka:87} for diffusions using the Malliavin
calculus, extended by Kohatsu-Higa~\cite{Kohatsu:03} for general
random variables on Wiener space, and adapted by
Bally~\cite{Bally:03} to deal with local ellipticity condition, do
not apply to our situation. Nevertheless, our approach follows the
same idea which consists in discretizing the time-interval and
writing the increments of the process on each subdivision interval
as the sum of a Gaussian term plus a remaining term. The
non-degeneracy of the Malliavin matrix proved by Gu\'erin implies
a deterministic lower bound for the smallest eigenvalue of the
Gaussian term covariance matrix. On the other hand, the upper
bound of the upper eigenvalue is random, due to the unboundedness
of the coefficients, and depends on the process itself, which
considerably complicates the problem. These estimates  on the
eigenvalues allow us to obtain a lower bound for the density of
the Gaussian term. In order to estimate the remaining term, we
need to refine some result on conditional Malliavin calculus to
deal with our specific situation. These results and our method
could be applied in other cases where the (invertible) Malliavin
covariance matrix of some functional has randomly upper-bounded
eigenvalues. The lower bound we finally obtain implies the
positivity of the solution of the Landau equation, and partially
answers to an analytical conjecture.

For the proof of the upper bound, we use tools of usual Malliavin
calculus. As the coefficients are not bounded, the proof differs
from the standard way to obtain Gaussian-type upper bounds. In
order to deal with a bounded  martingale term quadratic variation,
we consider the stochastic differential equation satisfied by some
logarithmic functional of the process. We then use an exponential
inequality for this martingale term. The diffusion matrix being
degenerate, we cannot apply Girsanov's theorem, which yields to
some unusual estimate.

The paper is organized as follows. In Section 2, we introduce the
Landau process as well as the main result. The relations with the
Fokker-Planck-Landau equation are also explained, as the
analytical interpretation of our results. In Section 3, we prove
general results on conditional Malliavin calculus. The proof of
the lower bound is given in Section $4$.
We finally show in Section 5 an upper-bound for the density. \vskip
12pt

In all the paper,  $C$ will denote an arbitrary constant whose
value may change from line to line.

\Section{The nonlinear Landau process and the main results}
\subsection{The nonlinear Landau process}

The Landau process is defined on  a filtered probability space
$(\Omega,{\mathcal F},({\mathcal F}_t)_{t\geq 0}, \mathbb{P})$. Fix $T>0$. We
consider $d\:$ independent space-time white noises
$W=(W^1,...,W^d)$  on $[0,T] \times [0,1]$, defined on
$(\Omega,{\mathcal F},\mathbb{P})$ and with covariance measure
$d\alpha dt$ on $[0,1]\times \mathbb{R}_+$ (cf. Walsh
\cite{Walsh:84}). Let $X_0$ be a random vector on $\mathbb{R}^d$,
independent of $W$. We denote by $(\mathcal{F}_t) _{t\geq 0}$ the
filtration generated by $W$ and $X_0$. In order to model the
nonlinearity, we also consider the probability space
$([0,1],\mathcal{B}([0,1]) ,d\alpha)$, $d\alpha$ denoting
Lebesgue measure. We denote by $\mathbb{E}$, $\mathbb{E}_{\alpha
}$ the expectations and $\mathcal{L}$, $\mathcal{L}_{\alpha }$ the
distributions of a random variable on
$(\Omega,\mathcal{F},\mathbb{P}) $, respectively on $([ 0,1]
,\mathcal{B}([0,1]),d\alpha)$. \vskip 12pt Let us consider the
following nonlinear stochastic differential equation.
\begin{defi}
A couple of processes $(X,Y) $ on $( \Omega ,\mathcal{F},(
\mathcal{F} _{t})_{t\geq 0},\mathbb{P}) \times ([0,1] ,
\mathcal{B}([ 0,1]) ,d\alpha)$ is defined as a solution of the
Landau stochastic differential equation if $\mathcal{L}\left(
X\right) =\mathcal{L}_{\alpha }\left( Y\right)$ and for any $t
\geq 0$,
\begin{equation}
\label{eq.landau} X_{t}=X_{0}+\int_{0}^{t}\int_{0}^{1}\sigma
\left( X_{s}-Y_{s}\left( \alpha \right) \right) \cdot W\left( d\alpha
,ds\right) +\int_{0}^{t}\int_{0}^{1}b\left( X_{s}-Y_{s}\left(
\alpha \right) \right) d\alpha ds,
\end{equation}
where $\sigma$ and $b$ are the coefficients of the spatially
homogeneous Landau equation for a generalization of Maxwellian
molecules (cf. Villani \textnormal{\cite{Villani2:98}}, Gu\'erin
\textnormal{\cite{Guerin:03}}).
\end{defi}
 More specifically, $\sigma$ is a $d \times d$
matrix (and $\sigma ^{\ast }$ denotes its adjoint matrix) such
that
\begin{equation*}
\sigma  \sigma ^{\ast }=a
\end{equation*}
where  $a$ is the $d \times d$ non-negative symmetric matrix given
by
\begin{equation}
\label{a} a_{ij}(z) =h(\vert z \vert^{2})(\vert z
\vert^{2}\delta_{ij}-z_{i}z_{j}),\qquad \forall\
(i,j)\in\{1,...,d\}^2
\end{equation}
($\delta_{ij}$ denotes the Kronecker symbol). Moreover,
\begin{equation*}
b_{i}(z) =\sum_{j=1}^{d} \partial_{z_j} a_{ij}(z) =-( d-1) h(\vert
z \vert^{2}) z_{i}, \qquad \forall\ i\in\{1,...,d\}.
\end{equation*}
When $h$ is a constant function, we recognize the coefficients of
the spatially homogeneous Landau equation for Maxwellian
molecules, cf.~\cite{Villani2:98}.

\vskip 12pt In all what follows, we assume the following
hypotheses: \vskip 12pt (H1): {\it The initial random variable
$X_0$ has finite  moments of order $k\geq 2$. } \vskip 12pt (H2):
{\it The function $h$ is defined on $\mathbb{R}_{+}$, sufficiently
smooth in order to get $\sigma$ and $b$ of class
$\mathcal{C}^{\infty}$ with bounded derivatives, and there exist
$m,M>0$ such that for all $r\in \mathbb{R}_{+}$,}
\begin{equation}
\label{boundh} m\leq h(r) \leq M.
\end{equation}

For example, in dimension two,
\begin{equation*}
\sigma \left( z\right) =\sqrt{h(|z|^{2}) }\left(
\begin{array}{cc}
z_{2}&0 \\
-z_{1}&0
\end{array}
\right),
\end{equation*}
 and in dimension three,
\begin{equation*}
\sigma \left( z\right) =\sqrt{h(|z| ^{2}) }\left(
\begin{array}{ccc}
z_{2} & -z_{3} & 0 \\
-z_{1} & 0 & z_{3} \\
0 & z_{1} & -z_{2}
\end{array}
\right),  \label{C2Strois}
\end{equation*}
and (H2) is satisfied for convenient function $h$.

\begin{defi}The $d$-dimensional stochastic process $X=(X_t, t \geq 0)$ is
called a nonlinear Landau process if there exists a process $Y$
defined on $[0,1]$ such that $(X,Y)$ is solution of the Landau SDE
(\ref{eq.landau}).
\end{defi}
This process has been introduced by Gu\'erin \cite{Guerin:02}
and~\cite{Guerin:03}, and gives a probabilistic interpretation of
the  spatially homogeneous Landau equation for generalized
Maxwellian molecules in the following sense.
\begin{prop}
If $(X,Y)$ is a solution of the Landau SDE
\textnormal{(\ref{eq.landau})}, then the family of laws $\left(
P_{t}\right) _{t\geq 0}$ of $\:(X_t)_{t\geq 0}$ (or of $\:(Y_t)
_{t\geq 0}$) satisfies for any $\varphi \in
\mathcal{C}_{b}^{2}\left( \mathbb{R}^{d}, \mathbb{R}\right)$,
\begin{eqnarray}
\frac{d}{dt} \int_{\mathbb{R}^d} \varphi \left( v\right)
P_{t}\left( dv\right) &=&\frac{1}{2}
\sum_{i,j=1}^{d}\int_{\mathbb{R}^{d}} \left( \int_{
\mathbb{R}^{d}} a_{ij}\left( v-v_{\ast }\right)P_{t}\left(
dv_{\ast }\right)
\right) \partial _{ij}\varphi \left( v\right)P_{t}\left( dv\right) \notag \\
&&+\sum_{i=1}^{d}\int_{\mathbb{R}^{d}} \left( \int_{
\mathbb{R}^{d}} b_{i}\left( v-v_{\ast }\right) P_{t}\left(
dv_{\ast }\right)\right) \partial _{i}\varphi \left(
v\right)P_{t}\left( dv\right). \label{C2landau-proba}
\end{eqnarray}
\end{prop}
 The proof is obtained using It\^o's Formula.
 \vskip 12pt

The equation (\ref{C2landau-proba}) is a weak form
 of the nonlinear partial differential equation
\begin{equation}
\frac{\partial f}{\partial t}\left( t,v\right) =\frac{1}{2}\sum_{i,j=1}^{d}%
\frac{\partial }{\partial v_{i}}\left\{
\int_{\mathbb{R}^{d}}a_{ij}\left( v-v_{\ast }\right) \left[ f\left(t, v_{\ast }\right) \frac{%
\partial f}{\partial v_{j}}\left( t,v\right) -f\left( t,v\right) \frac{%
\partial f}{\partial v_{\ast j}}\left(t, v_{\ast }\right) \right]dv_{\ast
} \right\}. \label{C2landau-initial}
\end{equation}
This equation is a spatially homogeneous Fokker-Planck-Landau
equation and models collisions of particles in a plasma. It can
also be obtained as limit of Boltzmann equations when collisions
become grazing
(\cite{Goudon:97},~\cite{Villani:98},~\cite{Guerin:04}). The
function $f(t,v) \geq 0$ is the density of particles with
velocity $v\in \mathbb{R}^{d}$ at time $t \geq 0$.

 \vskip 12pt The results proved by Gu\'erin
\cite{Guerin:02} can be summarized as follows.
\begin{thm}
\label{C2ExistenceProba} Fix $T>0$. Assume \textnormal{(H1)},
\textnormal{(H2)} and that the law of $X_{0}$ is not a Dirac
measure. Then there exists a unique couple $(X,Y)$ such that for any $p>1$,
$\mathbb{E}[\sup_{t\leq T}|X_t|^p]<+\infty$, solution of the Landau
SDE (\ref{eq.landau}).

Moreover, for any $t>0$, the regular version of the conditional
distribution of $X_{t}$ given $X_0$ is absolutely continuous with
respect to Lebesgue measure and  its density function
$f_{X_0}(t,v)$ is ($P_0$-a.s.) of class $\mathcal{C}^{\infty}$.
\end{thm}

For the proof of the existence and regularity of a density for
each $P_t, t>0$, Gu\'erin uses tools of Malliavin calculus, the
degeneracy of the matrix $\sigma$ being compensated by the effect
of  the nonlinearity.

Gu\'erin's result leads, using the probabilistic interpretation,
to the existence and uniqueness of a smooth solution for the
Landau equation, given by
$f(t,v)=\int_{\mathbb{R}^d}f_{x_0}(t,v)P_0(dx_0)$.

\subsection{The main results}

The aim of this paper is to obtain  some  upper and lower bounds
for the conditional density $f_{X_0}(t,v)$ of $X_t$ given $X_0$,
for any time $t$ in a bounded interval $(0,T]$. We deduce from
them the strict positivity of the density and some bounds and
positivity for the solution of the Landau equation. The research
of a lower bound for this equation was partially developed in
Villani \cite{Villani:98}. In that paper, the author obtained (in
Section 7-Theorem 3) a result in the case of Maxwellian molecules,
assuming that the initial condition is bounded below by a
Maxwellian function. The general case is much more complicated and
a conjecture was stated in \cite[Proposition 6]{Villani:98}, but
never proved. \vskip 12pt We now assume the additional
non-degeneracy hypothesis. \vskip 12pt (H3): {\it  For all $\xi
\in \mathbb{R}^d$, $\mathbb{E}[|X_0|^2 |\xi|^2-<X_0,\xi>^2]>0$.}
\vskip 12pt
\begin{rem}
\label{example} Hypothesis \textnormal{(H3)} means that the
support of the law of $X_0$ is not embedded in a line. In
particular, it holds for the two extreme cases, if either the law
$P_0$ of $X_0$ has a density $f_0$ with respect to Lebesgue
measure, or if $P_0=\frac{\delta_{x_1}+\delta_{x_2}}{2}$, with
$x_1$ and $x_2$ non collinear vectors.
\end{rem}

The main theorem of this article is the following :
\begin{thm} \label{maintheorem}
Fix $T>0$ and assume \textnormal{(H1)}, \textnormal{(H2)}.
\begin{itemize}
 \item[\textnormal{(a)}] Assume moreover \textnormal{(H3)}. Then for any $0<t \leq T$ and $v \in
\mathbb{R}^d$, there exist two constants $c_1(T,v,X_0)$ and
$c_2(T,v,X_0)$ (explicitely given in the proof), such that
$P_0$-a.s.,
$$
f_{X_0}(t,v) \geq c_1(T,v,X_0)\, t^{-d/2} e^{-c_2(T,v,X_0)
\frac{\vert v-X_0 \vert^2}{t}} .
$$
\item[\textnormal{(b)}] For any $0<t \leq T$ and $v \in
\mathbb{R}^d$, there exist constants $c_1(T),c_2(T),c_3(T,X_0)$
such that $P_0$-a.s.,
\begin{equation*}
f_{X_0}(t,v) \leq c_3(T,X_0) \, t^{-d/2} e^{-\frac{( \ln(1+
|v|^2)-\ln(1+ |X_0|^2)-c_1 t )^2}{c_2 t}}.
\end{equation*}
\end{itemize}
\end{thm}
\begin{cor}
For any $t>0$, the density function $f_{X_0}(t,v)$ is positive.
\end{cor}

 As a consequence of Theorem~\ref{maintheorem} and
writing $f(t,v)=\int_{\mathbb{R}^d} f_{x_0}(t,v)P_0(dx_0)$, we
obtain the positivity and bounds for the solution of the Landau
equation (\ref{C2landau-initial}). \vskip 12pt

We obtain (a) by adapting the approach of Kohatsu-Higa
\cite{Kohatsu:03}, in which a key tool is  conditioned Malliavin
calculus for general random processes with ellipticity and bounded
coefficients.  To deal with our degenerate process, we need
refined conditional Malliavin calculus, that will be given in the
next section.

\section{Conditional Malliavin calculus}

Recall some basic notions of the Malliavin calculus related to the
space-time white noise $W$.
 Fix $T>0$.  Let the Hilbert space
$\mathcal{H}={{L}}^2([0,T] \times [0,1]; \mathbb{R}^d)$. For any
$h \in \mathcal{H}$, we set $$W(h)= \int_0^T \int_0^1 h(r,z) \cdot
W(dr,dz).$$ Let $\mathcal{S}$ denote the class of smooth random
variables $F=f(W(h_1),...,W(h_n))$, where $h_1,...,h_n$ are in
$\mathcal{H}$, $n \geq 1$, and $f$ is of class
$\mathcal{C}^{\infty}$ on $\mathbb{R}^n$ with polynomial growth
derivatives.

Given $F$ in $S$, its derivative is the $d$-dimensional
stochastic process $DF=(D_{(r,z)} F=(D_{(r,z)}^{1}
F,...,D_{(r,z)}^{d} F), (r,z) \in [0,T] \times [0,1])$, where the
$D_{(r,z)} F$ are $\mathcal{H}$-valued random vectors given, for
$l=1,...,d$, by
\begin{equation*}
D_{(r,z)}^{l}F=\sum_{i=1}^{n} \partial_{x_i} f ( W(
h_1),...,W(h_n)) h^{l}_i(r,z).
\end{equation*}
More generally, if $F$ is a smooth random variable and $k$ is an
integer, set $D_{\alpha}^{(k)} F=D_{\alpha_1} \cdots D_{\alpha_k}
F$, where $\alpha=(\alpha_1,...,\alpha_k)$,
$\alpha_i=(r_i,z_i)\in[0,T] \times[0,1]$, for the $k$-th order
derivative of $F$. Then for every $p\geq1$ and any natural number
$m$, we denote by $\mathbb{D}^{m,p}$  the closure of $\mathcal{S}$
with respect to the semi-norm $\Vert \cdot \Vert_{m,p}$ defined by
\begin{equation*}
\Vert F \Vert _{m,p}=\biggl( \mathbb{E}[|F|^{p}] +\sum_{k=1}^{m}
\mathbb{E}[ \| D^{(k)} F \| _{\mathcal{H}^{\otimes
k}}^{p}]\biggr)^{1/p},
\end{equation*}
where
\begin{equation*}
\| D^{(k)}F \|_{\mathcal{H}^{\otimes k}}^2=\sum_{l_1,...,l_k=1}^d
\int \cdots \int_{([0,T] \times [0,1])^{k}} | D_{\alpha_1}^{l_{1}}
\cdots D_{\alpha_k}^{l_k} F | ^{2} \, d\alpha_1 \cdots d \alpha_k.
\end{equation*}

For any fixed $s \in [0,T]$, we define the conditional versions of
the Sobolev norms related to $W$  with respect to
${\mathcal{F}}_s$. Let  $p \geq 1$, and $n \geq 1$, $m \geq 0$
 natural integers. For any function $f \in
{{L}}^2(([0,T]\times [0,1])^n; \mathbb{R}^d)$ and any random
variable $F \in \mathbb{D}^{m,p}$, we define
\begin{align*}
\mathcal{H}_s&=L^2 ([s,T]\times
[0,1]; \mathbb{R}^d), \\
\Vert f \Vert_{\mathcal{H}_{s}^{\otimes n}}&= \biggl(
\int_{([s,T]\times [0,1])^n}
\vert f(r,z) \vert^2 dz_1 \cdots dz_n dr_1 \cdots dr_n \biggr)^{1/2}, \\
\Vert F \Vert_{m,p,s} &= \biggl( \mathbb{E}[|F|^p|{\mathcal{F}}_s]
+ \sum_{k=1}^m \mathbb{E}[\Vert D^{(k)} F
\Vert^p_{\mathcal{H}_{s}^{\otimes k}}|{\mathcal{F}}_s]
\biggr)^{1/p}.
\end{align*}
Moreover, we write $\gamma_{F}(s)$ for the Malliavin covariance
matrix with respect to $\mathcal{H}_s$, that is,
\begin{equation*}
\gamma_{F}(s)= ( \langle DF^i, D F^j
\rangle_{\mathcal{H}_s})_{1\leq i,j\leq d}.
\end{equation*}
For any  $u\in L^2(\Omega;\mathcal H)$ such that $u(r,z) \in
\mathbb{D}^{m,p}$, for all $(r,z) \in [0,T] \in [0,1]$, we define
\begin{equation*}
\Vert u \Vert_{m,p,s} = \biggl( \mathbb{E}[\Vert u
\Vert_{\mathcal{H}_{s}}^p|{\mathcal{F}}_s] + \sum_{k=1}^m
\mathbb{E}[\Vert D^{(k)} u \Vert ^p_{\mathcal{H}_{s}^{\otimes
k+1}}|{\mathcal{F}}_s] \biggr)^{1/p}.
\end{equation*}
 We denote by $\delta$
the adjoint of the operator $D$, which is an unbounded operator on
${{L}}^2(\Omega; \mathcal{H})$ taking values in ${{L}}^2(\Omega)$
(see~\cite[Def.$1.3.1$]{Nualart:95}). In particular, if $u$
belongs to Dom~$\delta$, then $\delta (u)$ is the element of
${{L}}^2 (\Omega)$ characterized by the following duality
relation:
\begin{equation*}
\mathbb{E}[F \delta (u)]= \mathbb{E} \biggl[\int_0^T \int_0^1
D_{(r,z)} F \cdot u(r,z) dz dr \biggl], \; \; \text{for any} \; F
\in \mathbb{D}^{1,2}.
\end{equation*}

With this notation one has the following estimate for the
conditional norm of the operator $\delta$ (cf.
{\cite[(2.15)]{Moret:01}}):
\begin{equation} \label{3.10}
\Vert \delta(u{\bf 1}_{[s,T]\times [0,1]})\Vert_{m,p,s} \leq
c_{m,p} \Vert u \Vert_{m+1,p,s},
\end{equation}
for some constant $c_{m,p}>0$.

\vskip 12pt  We next give a conditional version of the integration
by parts formula. The proof follows similarly as the
non-conditional version (cf.~\cite[Proposition 3.2.1]{Nualart:98},
and is therefore omitted.
\begin{prop} \label{ipfc}
Fix $n\geq1$. Let $F, Z_s, G \in (\cap_{p\geq 1}\cap_{m\geq 0}
\mathbb{D}^{m,p})^d$ be three random vectors where $Z_s$ is
${\mathcal{F}}_s$-measurable and such that $(\textnormal{det} \,
\gamma_{F+Z_s}(s))^{-1}$ has finite moments of all orders. Let $g
\in \mathcal{C}^{\infty}_p(\mathbb{R}^d)$. Then, for any
multi-index $\alpha=(\alpha_1,...,\alpha_n) \in \{1,\dots,d \}^n$,
there exists an element $H^s_{\alpha}(F,G) \in \cap_{p\geq
1}\cap_{m\geq 0} \mathbb{D}^{m,p}$ such that
\begin{equation*}
\mathbb{E}[ (\partial_{\alpha}g) (F+Z_s) G|{\mathcal{F}}_s]=
\mathbb{E}[ g(F+Z_s) H^s_{\alpha}(F,G)|{\mathcal{F}}_s],
\end{equation*}
where the random variables $H^s_{\alpha}(F,G)$ are recursively
given by
\begin{align*}
H^s_{(i)}(F,G)& =\sum_{j=1}^d\delta(G \, (\gamma_{F}(s)^{-1})_{i j} \, DF^j), \\
H^s_{\alpha}(F,G)& =
H^s_{(\alpha_n)}(F,H^s_{(\alpha_1,\dots,\alpha_{n-1})}(F,G)).
\end{align*}
\end{prop}

As a consequence of this integration by parts formula, one derives
the following expression for the conditional density given
${\mathcal{F}}_s$ of a random vector on the Wiener space, in a
similar way as in \cite[Proposition 4]{Moret:01}.
\begin{cor} \label{expression}
Let $F\in (\cap_{p\geq 1}\cap_{m\geq 0} \mathbb{D}^{m,p})^d$ be a
random vector such that $(\textnormal{det} \, \gamma_{F}(s))^{-1}$
has finite moments of all orders. Let $P_s$ and $p_s$ denote,
respectively, the conditional distribution and density of $F$
given ${\mathcal{F}}_s$. Let $\sigma$ be a subset of the set of
indices of $\{1,...,d\}$. Then, for any $v \in \mathbb{R}^d$,
$P_s$-a.s.
$$
p_s(v)=(-1)^{d-\vert \sigma \vert} \mathbb{E} [1_{\{F^i>v_i,\,
i\in \sigma\,; F^i<v_i,\, i\notin \sigma\,;\, i=1,...,d\}}
 H^s_{(1,...,d)}(F,1) |{\mathcal{F}}_s],
$$
where $\vert \sigma \vert$ denotes the cardinality of $\sigma$.
\end{cor}

The next result gives a precise estimate of the Sobolev norm of
the random variables $H^s_{\alpha}(F,G)$.
\begin{prop} \label{normHc}
 Let $F \in (\cap_{p \geq 1} \cap_{m\geq 0}
\mathbb{D}^{m,p})^d$ and $G \in \cap_{p \geq 1} \cap_{m\geq 0}
\mathbb{D}^{m,p}$ be two  random vectors such that
$(\textnormal{det} \, \gamma_F(s))^{-1}$ has finite moments of all
orders. Assume that there exist positive
$\mathcal{F}_s$-measurable finite random variables $Z_s$ and $Y_s$
(eventually deterministic) such that for all $p>1$ and $m\geq 1$,
\begin{align} \label{a1}
\mathbb{E}[\Vert D^{(m)}(F^i) \Vert^p_{\mathcal{H}_s^{\otimes m}}
|{\mathcal{F}}_s]^{1/p} &\leq c_1(m,p) Z_s, \; \; i=1,...,d;
 \\
\label{a2} \mathbb{E}[(\textnormal{det} \,\gamma_{F}(s)
)^{-p}|{\mathcal{F}}_s]^{1/p} &\leq c_2(p) Z_s^{-2d} Y_s,
\end{align}
where $c_1(m,p)$ and $c_2(p)$ are positive constants. Then, for
any multi-index $\alpha=(\alpha_1,...,\alpha_n) \in \{1,\dots,d
\}^n$, $n \geq 1$, there exists a constant $C>0$ (depending on
$m$, $p$, $\alpha$, $T$), such that
\begin{equation*}
\Vert H^s_{\alpha}(F,G) \Vert_{0,2,s} \leq C \Vert G
\Vert_{n,2^{n+1},s} Z_s^{-n} \prod_{i=1}^n \bigg(\sum_{j=1}^{i+1}
(Y_s)^j\bigg).
\end{equation*}
\end{prop}

\begin{proof}
The proof of this result follows  the iteration argument appearing
in the proof of~\cite[Lemma 12]{Moret:01} or~\cite[Lemma
4.11]{Dalang:03}, but in a general setting. That is, we use
(\ref{3.10}) and H\"older's inequality for the conditional
Malliavin norms (cf.~\cite[Proposition 1.10, p.50]{Watanabe:84} to
obtain
\begin{align}\label{induction}\nonumber
\Vert H^s_{\alpha}(F,G) \Vert_{0,2,s}
&=\Vert \sum_{j=1}^d \delta(H^s_{(\alpha_1,...,\alpha_{n-1})}(F,G) \,
(\gamma_{F}(s)^{-1})_{\alpha_n j} \, DF^j) \Vert_{0,2,s} \\
&\leq C \Vert H^s_{(\alpha_1,...,\alpha_{n-1})}(F,G)
\Vert_{1,2^2,s} \sum_{j=1}^d \Vert (\gamma_{F}(s)^{-1})_{\alpha_n
j} \Vert_{1,2^3,s} \, \Vert D(F^j) \Vert_{1,2^3,s}.
\end{align}

Note that, as proved in \cite[Lemma 11]{Dalang:03}, for $m\geq 1$ and $p>1$,
\begin{eqnarray} \label{a6} \nonumber
&&\mathbb{E} [\Vert D^{(m)} (\gamma_F(s))_{ij}
\Vert^{p}_{\mathcal{H}^{\otimes m}_{s}}|{\mathcal{F}}_s]
 = \mathbb{E} [\Vert D^{(m)} (  \langle D(F^i),
 D(F^j) \rangle_{\mathcal{H}_s} ) \Vert^{p}_{\mathcal{H}^{\otimes m}_{s}}|{\mathcal{F}}_s] \\
\nonumber && \leq  C \sum_{l=0}^m \left(\begin{array}{c} \!\!m\!\!
\\  \!\!l\! \! \end{array} \right)^p
\{ (\mathbb{E}[ \Vert D^{(l+1)}(F^i)
\Vert^{2p}_{\mathcal{H}^{\otimes (l+1)}_{s}} |
{\mathcal{F}}_s ])^{1/2} \\
&&\qquad \qquad \qquad  \qquad \qquad \qquad\times
(\mathbb{E}[\Vert D^{(m-l+1)} (F^j)
\Vert^{2p}_{\mathcal{H}^{\otimes (m-l+1)}_{s}} |{\mathcal{F}}_s
])^{1/2} \}.\nonumber
\end{eqnarray}

Therefore, by (\ref{a1}) we get,  for $1\leq i,j\leq d$, \be
\label{dergamma}
 \|D((\gamma_{F}(s))_{ij})\|_{m,p,s}\leq C Z_s^2.
 \ee

Now, Cramer's formula gives
$$
|(\gamma_{F}(s)^{-1})_{ij}|=|A_{ij}(s)(\text{det} \,
\gamma_F(s))^{-1}|,
$$
where $A_{ij}(s)$ denotes the cofactor of $(\gamma_{F}(s))_{ij}$.
By some straightforward computations, it is easily checked that
there exists a constant $C>0$ such that
$$
\vert A_{ij}(s) \vert \leq C \Vert D(F)
\Vert_{\mathcal{H}_s}^{2(d-1)}.
$$
Therefore, Cauchy-Schwarz inequality for conditional expectations
and hypotheses (\ref{a1}) and (\ref{a2}) yield
\begin{align} \label{a3} \nonumber
(\mathbb{E}[ ((\gamma_{F}(s)^{-1})_{ij})^p|{\mathcal{F}}_s ])^{1/p} &\leq
C (\mathbb{E}[ \Vert D(F)
 \Vert_{\mathcal{H}_s}^{4p (d-1)}|{\mathcal{F}}_s])^{1/(2p)} \times (\mathbb{E}[ (\text{det} \,
\gamma_F(s))^{-2p}|{\mathcal{F}}_s])^{1/(2p)} \\
 & \leq C Z_s^{2(d-1)} Z_2^{-2d} Y_s =C Z_s^{-2} Y_s.
\end{align}

Iterating the equality
\begin{equation*}
D ( \gamma_{F}(s)^{-1} )_{ij}=- \sum_{k,l=1}^d ( \gamma_{F}(s)^{-1}
)_{ik} D ( \gamma_{F}(s))_{kl} ( \gamma_{F}(s)^{-1} )_{jl},
\end{equation*}
and using H\"older's inequality for conditional expectations, we
obtain
\begin{eqnarray} \label{a4} \nonumber
&& \sup_{i,j} \mathbb{E} [
\Vert D^{(m)} ((\gamma_{F}(s))^{-1} )_{ij} \Vert^p_{\mathcal{H}_s^{\otimes m}}|
{\mathcal{F}}_s] \\
\nonumber &&\leq  C \sup \sum_{r=1}^m \sum_{\scriptstyle m_1+
\cdots +m_r=m \atop \scriptstyle m_l \geq 1, \, l=1,...,r}
\mathbb{E} [ \Vert D^{(m_1)} (\gamma_F (s))_{i_1j_1}
\Vert^{p(r+1)}_{\mathcal{H}_{s}^{\otimes
m_1}}|{\mathcal{F}}_s]^{1/(r+1)}
 \times \cdots \\ \nonumber
&&\qquad \qquad \qquad  \qquad\qquad \times \mathbb{E} [\Vert
D^{(m_r)} (\gamma_F (s))_{i_rj_r}
\Vert^{p(1+r)}_{\mathcal{H}_{s}^{\otimes m_r}}|{\mathcal{F}}_s]^{1/(r+1)} \\
&& \qquad \qquad \qquad \qquad \qquad\times
\sup_{i,j} \mathbb{E} [| ((\gamma_F (s))^{-1} )_{ij}|^{p (r+1)^2}|{\mathcal{F}}_s]^{1/(r+1)},
\end{eqnarray}
where the supremum before the summation is over
$i_1,j_1,...,i_{2r+1},j_{2r+1} \in \{ 1,...,d\}$.

Introducing (\ref{dergamma}) and (\ref{a3}) into (\ref{a4}) gives
\begin{equation} \label{a5}
\Vert D (\gamma_{F}(s)^{-1})_{ij} \Vert_{m,p,s} \leq C
Z_s^{-2}\sum_{r=1}^m Y_s^{r+1}.
\end{equation}
and thus
$$
\Vert (\gamma_{F}(s)^{-1})_{ij} \Vert_{m,p,s} \leq C
Z_s^{-2}\sum_{r=0}^m Y_s^{r+1}.
$$
Therefore, iterating $n$ times (\ref{induction}),
  it yields
\begin{align*}
\Vert H^s_{\alpha}(F,G) \Vert_{0,2,s} &\leq
C \Vert H^s_{(\alpha_1,...,\alpha_{n-1})}(F,G) \Vert_{1,2^2,s} Z_s^{-1} (Y_s+Y_s^2) \\
& \leq C \Vert H^s_{(\alpha_1)}(F,G) \Vert_{n-1,2^n,s} Z_s^{-n+1}
\prod_{i=1}^{n-1}
(\sum_{j=1}^{i+1} Y_s^j)\\
& \leq C \Vert G \Vert_{n,2^{n+1},s} Z_s^{-n} \prod_{i=1}^n
(\sum_{j=1}^{i+1} Y_s^j),
\end{align*}
which concludes the proof of the Proposition.
\end{proof}

\vskip 12pt The last result of this section will be used later in
order to prove condition (\ref{a1}) of Proposition~\ref{normHc}
when $F$ is the Landau random variable $X_t$.
\begin{prop} \label{deter}
Fix $\epsilon_0>0$ and $0<\alpha_1<\alpha_2$. Fix $c_1>0$ and for
$q>1$, let $c_2(q)$ be finite. Let $Z$ be a positive random
variable such that for all $\epsilon \leq \epsilon_0$, there exist
two random variables $X(\epsilon)$, $Y(\epsilon)$ such that $Z
\geq X(\epsilon) -Y(\epsilon)$ a.s., and
\begin{enumerate}
\item[$(1)$] $X(\epsilon) \geq c_1 \, \epsilon^{\alpha_1}$ \; a.s., and
\item[$(2)$] there exists a positive $\mathcal{F}_s$-measurable finite random
variable $G_s$ (eventually deterministic) such that for any $q>1$,
$\mathbb{E}[\vert Y(\epsilon) \vert^q|{\mathcal{F}}_s] \leq  c_2(q) \, \epsilon^{q \, \alpha_2} G_s^q$.
\end{enumerate}
Then, for any $p\geq 1$ and $q>\frac{p \alpha_1}{\alpha_2-\alpha_1}$, there exists a constant $c_3$ depending on $c_1, c_2(q), \alpha_1, \alpha_2$, but not on $Z$, $G_s$ or $\epsilon_0$ such that, a.s.,
$$
\mathbb{E}[Z^{-p}|{\mathcal{F}}_s] \leq c_3 \, \epsilon_0^{-p \, \alpha_1}(1+ \epsilon_0^{q(\alpha_2-\alpha_1)} G_s^q).
$$
\end{prop}

\begin{proof}
For $p\geq1$, we write
\begin{equation} \label{ol}
\mathbb{E}  [Z^{-p}|{\mathcal{F}}_s]=\int_0^{\infty} p y^{p-1} \mathbb{P} \{
Z^{-1}>y |{\mathcal{F}}_s\} \, dy.
\end{equation}
Let $k=(\frac{c_1}{2} {\epsilon_0}^{\alpha_1})^{-1}$. For $y \geq k$, let $\epsilon=(\frac{2}{c_1})^{1/\alpha_1} y^{-1/\alpha_1}$.
Then $\epsilon \leq \epsilon_0$ and $y^{-1}=\frac{c_1}{2} \epsilon^{\alpha_1}$. By Chebychev's inequality with $q>1$,
\begin{align*}
\mathbb{P} \{ Z^{-1}>y |{\mathcal{F}}_s\} &
\leq  \mathbb{P} \{ Y_{\epsilon} > X_{\epsilon}- y^{-1} | {\mathcal{F}}_s\}
\leq \mathbb{P} \{Y_{\epsilon}  > \frac{c_1}{2} \, \epsilon^{\alpha_1} | {\mathcal{F}}_s \} \\
& \leq (\frac{c_1}{2} \, \epsilon^{\alpha_1})^{-q} \mathbb{E}[\vert Y_{\epsilon} \vert^q| {\mathcal{F}}_s] \\
&\leq (\frac{c_1}{2} \, \epsilon^{\alpha_1})^{-q} c_2(q) \epsilon^{q \, \alpha_2} G_s^q\\
&=c_q \epsilon^{q(\alpha_2-\alpha_1)} G_s^q=\tilde{c}_q y^{-q(\alpha_2-\alpha_1)/\alpha_1} G_s^q.
\end{align*}

Now, splitting the integral in
(\ref{ol}) into an integral over $[0,k]$ and another on $(k, +\infty)$, introducing this last inequality into (\ref{ol}) and choosing $q>\frac{p \alpha_1}{\alpha_2-\alpha_1}$, we obtain
\begin{align*}
\mathbb{E} [Z^{-p}|{\mathcal{F}}_s]&\leq k^{p} + p \int_k^{\infty} y^{p-1} \mathbb{P} \{
Z^{-1}>y |{\mathcal{F}}_s\} dy \\
&\leq c_p \epsilon_0^{-p \, \alpha_1} + c_{p,q} \int_k^{\infty} y^{p-1-q(\alpha_2-\alpha_1)/\alpha_1} G_s^q dy \\
&= c_p \epsilon_0^{-p \, \alpha_1} + c_{p,q} \epsilon_0^{-p \, \alpha_1+q(\alpha_2-\alpha_1)} G_s^q \\
&\leq c_3 \, \epsilon_0^{-p \, \alpha_1}(1+ \epsilon_0^{q(\alpha_2-\alpha_1)} G_s^q),
\end{align*}
which concludes the proof of the Proposition.
\end{proof}

\section{The Lower Bound}

 The aim of this section is to prove the lower
bound of Theorem \ref{maintheorem}. As in Kusuoka-Stroock
\cite{Kusuoka:87} and Kohatsu-Higa \cite{Kohatsu:03}, we
discretize the time interval $[0,t]$ and write $X_t$ as the sum of
a Gaussian term plus a remaining term. The lower bound for the
density of our process is deduced from a lower estimate of the
density of the Gaussian term and a technical part consists in the
choice of the discretization mesh in order to control the
remaining term. These steps can not be obtained from
\cite{Kusuoka:87} and \cite{Kohatsu:03}, as the eigenvalues of the
covariance matrix  of the Gaussian term are not bounded, but only
dominated by a random functional of the diffusion, due to the
unboundedness of the coefficients. We will then use the results on
conditional Malliavin calculus of the previous section.

\subsection{The Discretized  Process}
We want to obtain a lower bound  of the conditional density of the
Landau process with respect to the initial condition $X_0$, on
some finite interval $[0,T]$. Then, in all what follows, $X_0$
will be considered as a parameter, even if it is random, and all
the estimates we get will concern conditional expectations with
respect to this initial condition $X_0$.

Let $T>0$ and fix  $t\in (0,T]$. Let us introduce a natural
integer $N$, measurably depending on $X_0$, which will be chosen
later.

 Consider a time grid $0=t_0 < t_1 < \cdots < t_N=t$ and
let $\Delta=t_k-t_{k-1}=\frac{t}{N}$. We define the following
discretized sequence,
\begin{equation} \label{evolu}
X_{t_k}=X_{t_{k-1}}+J_k+\Gamma_k,
\end{equation}
where
$$ J_k= \int_{t_{k-1}}^{t_k} \int_0^1
 \sigma(X_{t_{k-1}}-Y_{t_{k-1}}(\alpha)) \cdot W(d\alpha, ds),
$$ and
\begin{align*}
 \Gamma_k&=\int_{t_{k-1}}^{t_k} \int_0^1
(\sigma(X_s-Y_s(\alpha))-
\sigma(X_{t_{k-1}}-Y_{t_{k-1}}(\alpha))) \cdot W(d\alpha, ds) \\
& \qquad +\int_{t_{k-1}}^{t_k} \int_0^1 b(X_s-Y_s(\alpha)) \,
d\alpha ds.
\end{align*}
Conditioned  with respect to $\mathcal{F}_{t_{k-1}}$, the random
variable $J_k$ is Gaussian with covariance matrix given by
$$
\Sigma(J_k)=(t_k-t_{k-1}) \int_0^1 a(X_{t_{k-1}}-Y_{t_{k-1}}(\alpha)) d\alpha.
$$

We wish to obtain a lower bound for the conditional density of the
random variable $X_{t_k}$ given $\mathcal{F}_{t_{k-1}}$. This will allow us to
prove the desired lower bound for the density of $X_t$ by a recursive method. Note that
from Theorem~\ref{C2ExistenceProba} this conditional density
exists and, from Watanabe's notation, can be written
$\mathbb{E}[\delta_z(X_{t_k}) |{\mathcal{F}_{t_{k-1}}} ]$, where
$\delta_z$ denotes the Dirac measure at the point $z \in
\mathbb{R}^d$.

We consider the following approximation of $\delta_z$. Let $\phi
\in \mathcal{C}^{\infty}_b(\mathbb{R}^d)$, $0\leq \phi \leq 1$,
$\int \phi=1$ and $\phi(x)=0$ for $\vert x \vert > 1$. For $\eta>0$, let
$$
\phi_{\eta}(x)= \eta^{-d} \phi(\eta^{-1} x).
$$
Remark that $\phi_{\eta}(x)=0$ for $\vert x \vert > \eta$. \vskip
12pt Our goal is to find a lower bound for the quantity
$\mathbb{E}[\phi_{\eta}(X_{t_{k}}-z) | {\mathcal{F}}_{t_{k-1}}]$,
independent of $\eta$. Let us apply the mean value theorem. We
have
\begin{align} \label{meanvalue} \nonumber
\mathbb{E}[\phi_{\eta}(X_{t_k}-z) | {\mathcal{F}}_{t_{k-1}}]& =
\mathbb{E}[\phi_{\eta} (X_{t_{k-1}}+J_k-z) |
{\mathcal{F}}_{t_{k-1}}] \\ \nonumber &  \qquad \qquad+
\sum_{i=1}^d \int_0^1 \mathbb{E} \biggl[\partial_{x^i} \phi_{\eta}
(X_{t_{k-1}}+J_k-z+\rho \Gamma_k) \Gamma^i_k |
{\mathcal{F}}_{t_{k-1}} \biggr] d \rho \\ \nonumber &
\geq\mathbb{E}[\phi_{\eta} (X_{t_{k-1}}+J_k-z) |
{\mathcal{F}}_{t_{k-1}}] \\ \nonumber &  \qquad \qquad- \bigg\vert
\sum_{i=1}^d \int_0^1 \mathbb{E} \biggl[\partial_{x^i} \phi_{\eta}
(X_{t_{k-1}}+J_k-z+\rho \Gamma_k) \Gamma^i_k |
{\mathcal{F}}_{t_{k-1}} \biggr] d \rho \bigg\vert \\
\end{align}
The two next subsections are devoted to obtain a lower bound for
the Gaussian term $\ \mathbb{E}[\phi_{\eta} (X_{t_{k-1}}+J_k-z) |
{\mathcal{F}}_{t_{k-1}}]\ $ and an upper bound for the remaining
term $$\vert \sum_{i=1}^d \int_0^1 \mathbb{E} [\partial_{x^i}
\phi_{\eta} (X_{t_{k-1}}+J_k-z+\rho \Gamma_k) \Gamma^i_k |
{\mathcal{F}}_{t_{k-1}}] d \rho \vert$$ of the RHS term of
(\ref{meanvalue}).

\subsection{Lower bound for the Gaussian term}

The following proposition gives a lower bound for the lower
eigenvalue and an upperbound for the upper eigenvalue of the matrix $\Sigma(J_k)$.

\begin{prop}  Under hypotheses \textnormal{(H1)}, \textnormal{(H2)}, \textnormal{(H3)},
 there exist two positive
constants $\lambda_1$ and $\lambda_2$ depending on $T$ such that
for any $k\in \{1,...,N\}$, almost surely,
\begin{eqnarray} \label{ev.min}
&&\inf_{\xi \in \mathbb{R}^d, \vert \xi \vert=1} \xi^{\ast}
\Sigma(J_k)
\xi \geq \lambda_1 \Delta \ ;\\
&& \label{ev.maj} \sup_{\xi \in \mathbb{R}^d, \vert \xi \vert =1}
\xi^{\ast} \Sigma(J_k) \xi \leq \lambda_2 \Delta (1+\vert
X_{t_{k-1}}\vert)^2.
\end{eqnarray}
 \end{prop}
\begin{proof}
In \cite{Guerin:02}, Gu\'erin shows that for each $\xi \in
\mathbb{R}^d$, one has
\begin{equation*}
\xi^{\ast}  \Sigma(J_k) \xi \geq \Delta m F(\xi,t_{k-1}),
\end{equation*}
where
\begin{equation*}
F(\xi,t)=\mathbb{E}[\vert X_t \vert^2 \vert \xi \vert^2- \langle
X_t, \xi \rangle^2],
\end{equation*}
and $m$ is defined in (\ref{boundh}).

 By Cauchy-Schwarz inequality, $F(\xi,t)$ is nonnegative, and
 since the law of $X_t$ has a density, $F(\xi,t)>0$ for any $t>0$ and $\xi\neq 0$.
Moreover, by Hypothesis (H3), this holds for $t\geq 0$. Then, as
the function $F(\xi,t)$ is positive and continuous on the compact
set $[0,T] \times \{\xi \in \mathbb{R}^d: \vert \xi \vert=1 \}$, a
strictly positive minimum is reached on this set.

Hence, for all $\xi \in \mathbb{R}^d$, $\vert \xi \vert=1$, we get
$$
\xi^{\ast}  \Sigma(J_k) \xi \geq \lambda_1 \Delta,
$$
where $\lambda_1>0$ is independent of $k$. That proves
(\ref{ev.min}).

Using the Lipschitz property of $\sigma$ (with Lipschitz constant
$C_{\sigma}$), we also obtain
\begin{eqnarray*}
\xi^{\ast} \Sigma(J_k) \xi&\leq& \Delta 2 C_{\sigma}^2 \int_0^1
(\vert
X_{t_{k-1}}\vert^2+\vert Y_{t_{k-1}}(\alpha) \vert^2) d\alpha \\
&=&\Delta 2 C_{\sigma}^2 (\vert
X_{t_{k-1}}\vert^2+\mathbb{E}[\vert X_{t_{k-1}}\vert^2])  \\
&\leq & \Delta 2 C_{\sigma}^2 (\vert
X_{t_{k-1}}\vert^2+\mathbb{E}[\sup_{0\leq s \leq T} \vert X_s\vert^2])  \\
&\leq& \lambda_2 \Delta (1+\vert X_{t_{k-1}}\vert)^2,
\end{eqnarray*}
and deduce (\ref{ev.maj}).
\end{proof}
\vskip 12pt The next result proves a lower bound for the
conditional density of the Gaussian term $X_{t_{k-1}}+J_k$ given
$\mathcal{F}_{t_{k-1}}$.
\begin{prop} \label{previouslowerbound}
Assume $0<\eta\leq \sqrt{\lambda_1 \Delta}$, and let $k\in \{1,...,N\}$. Then for  $(w,z)\in \Omega\times \mathbb{R}^d$
satisfying $ \vert X_{t_{k-1}}(\omega)-z \vert \leq
\sqrt{\lambda_1 \Delta}$, we get a.s.
$$
\mathbb{E}[\phi_{\eta}
(X_{t_{k-1}}+J_k-z)|{\mathcal{F}}_{t_{k-1}}]
 \geq \frac{1}{C_1 \Delta^{d/2} (1+\vert
 X_{t_{k-1}}\vert)^{d}},
$$
where $C_1:=e^2 (2 \pi)^{d/2} \lambda_2^{d/2}$.
\end{prop}

\begin{proof}
As $J_k$ is Gaussian,
\begin{align*}
& \mathbb{E}[\phi_{\eta} (X_{t_{k-1}}+J_k-z)|
{\mathcal{F}}_{t_{k-1}}] \\
&=\int_{\mathbb{R}^d} \phi_{\eta}(X_{t_{k-1}}+x-z) \frac{1}{(2
\pi)^{d/2} \text{det}(\Sigma(J_k))^{1/2}} \exp \biggl(-\frac{x^*
\Sigma(J_k)^{-1} x}{2}\biggr) \, dx \\
&=\int_{\mathbb{R}^d} \phi_{\eta}(\tilde{z}) \frac{1}{(2
\pi)^{d/2} \text{det}(\Sigma(J_k))^{1/2}}  \\
& \qquad \qquad \qquad \qquad \times\exp
\biggl(-\frac{(\tilde{z}+z-X_{t_{k-1}})^* \Sigma(J_k)^{-1}
(\tilde{z}+z-X_{t_{k-1}})}{2}\biggr) \, d\tilde{z}.
\end{align*}
Since $\vert \tilde{z} \vert \leq \eta\leq \sqrt{\lambda_1
\Delta}$, and using the assumption on $(\omega,z)$,
$$
\vert \tilde{z}+z-X_{t_{k-1}}\vert^2 \leq 2 \vert \tilde{z}
\vert^2 + 2\vert z - X_{t_{k-1}}\vert^2 \leq 4 \lambda_1 \Delta.
$$
Then, using (\ref{ev.min}) and (\ref{ev.maj}), we obtain
$$
\mathbb{E}[\phi_{\eta} (X_{t_{k-1}}+J_k-z)|
{\mathcal{F}}_{t_{k-1}}] \geq \frac{1}{C_1 \Delta^{d/2} (1+\vert
X_{t_{k-1}} \vert)^{d}},
$$
where $C_1:=e^2 (2 \pi)^{d/2} \lambda_2^{d/2}$.
\end{proof}

\subsection{Upper bound for the remaining term}

The key point consists in applying the conditional integration by
parts formula to the remaining term in (\ref{meanvalue}), taking
into account that $\int \phi =1$.  Then, in order to  obtain an
upper bound, we need to prove estimates for the conditional
Sobolev norms given ${\mathcal{F}}_{t_{k-1}}$ of the terms $J_k$
and $\Gamma_k$ of the discretized sequence (\ref{evolu}). Note
that as the coefficients of the Landau equation are unbounded,
these conditional bounds will depend on the random variable
$X_{t_{k-1}}$.

\begin{lem} \label{p1derivativegamma}
For any $p>1$, there exists a finite constant $C_T$ such that,
for $i\in \{1,...,d\}$ and $k\in \{1,...,N\}$,
$$
(\mathbb{E} [\vert\Gamma_k^i\vert^p| {\mathcal{F}}_{t_{k-1}}])^{1/p} \leq
C_T \Delta (1+\vert X_{t_{k-1}} \vert).
$$
\end{lem}

\begin{proof}
Note that $\, \mathbb{E} [\vert\Gamma_k^i\vert^p|
{\mathcal{F}}_{t_{k-1}}] \leq 2^{p-1}(A_1+A_2)$, where
\begin{align*}
A_1 &:=\mathbb{E} \biggl[ \biggl( \int_{t_{k-1}}^{t_k} \int_0^1
\sum_{j=1}^d (\sigma_{ij}(X_s-Y_s(\alpha))-
\sigma_{ij}(X_{t_{k-1}}-Y_{t_{k-1}}(\alpha))) \, W^j(d\alpha,
ds)\biggr)^p  |
{\mathcal{F}}_{t_{k-1}}\biggl], \\
A_2 &:= \mathbb{E}\biggl[ \biggl(\int_{t_{k-1}}^{t_k} \int_0^1
b_i(X_s-Y_s(\alpha)) d\alpha ds \biggr)^{p} |
{\mathcal{F}}_{t_{k-1}}\biggr].
\end{align*}
Using Burkholder's inequality for conditional expectations, we get
$$
A_1 \leq C \mathbb{E} \biggl[ \biggl( \int_{t_{k-1}}^{t_k}
\int_0^1 \sum_{j=1}^d (\sigma_{ij}(X_s-Y_s(\alpha))-
\sigma_{ij}(X_{t_{k-1}}-Y_{t_{k-1}}(\alpha)))^2 \, d\alpha ds
\biggr)^{p/2} | {\mathcal{F}}_{t_{k-1}}\biggl],
$$
and, from H\"older's inequality and the Lipschitz property of
$\sigma$, it yields
$$
A_1 \leq C \Delta^{p/2-1} \int_{t_{k-1}}^{t_k} \biggl(
\mathbb{E} [\vert X_s-X_{t_{k-1}} \vert^p |{\mathcal{F}}_{t_{k-1}}
]+ \mathbb{E} [\vert X_s -X_{t_{k-1}}\vert^p ]\biggr) ds.
$$
We now apply Burkholder's inequality and Lipschitz property, to obtain that, for $s \leq t_k$,
\begin{align*}
&\mathbb{E}[\vert X_s- X_{t_{k-1}}\vert^p|{\mathcal{F}}_{t_{k-1}}]
\leq C \Delta^{p/2-1} \biggl\{ \int_{t_{k-1}}^s \int_0^1
\mathbb{E}[\vert X_u
\vert^p + \vert Y_u(\alpha) \vert^p |{\mathcal{F}}_{t_{k-1}}] d \alpha du \\
&  \qquad \qquad \qquad \qquad \qquad \qquad + \Delta^{p/2}
\int_{t_{k-1}}^s \int_0^1 \mathbb{E}[\vert X_u \vert^p + \vert
Y_u(\alpha) \vert^p |{\mathcal{F}}_{t_{k-1}}] d
\alpha du  \biggr\} \\
&\qquad \qquad \qquad \leq C_T \Delta^{p/2-1}
\biggl(\int_{t_{k-1}}^s \mathbb{E}[\vert X_u \vert^p
|{\mathcal{F}}_{t_{k-1}}]+
\mathbb{E}[\vert X_u \vert^p] du \biggr) \\
&  \qquad \qquad \qquad \leq C_T \Delta^{p/2-1} \int_{t_{k-1}}^s
\mathbb{E}[\vert X_u -X_{t_{k-1}} \vert^p
|{\mathcal{F}}_{t_{k-1}}] du + C_T \Delta^{p/2}
(1+|X_{t_{k-1}}|)^p.
\end{align*}
By Gronwall's Lemma,
\begin{equation} \label{estimateA1}
\mathbb{E}[\vert X_s- X_{t_{k-1}} \vert^p |
{\mathcal{F}}_{t_{k-1}}] \leq C_T \Delta^{p/2} (1+\vert X_{t_{k-1}}
\vert)^p.
\end{equation}
Therefore,
\begin{equation} \label{estimateA}
A_1 \leq C_T \Delta^p (1+\vert X_{t_{k-1}} \vert)^p.
\end{equation}

On the other hand, using H\"older's inequality and Lipschitz property of
$b$, we have that
\begin{align*}
A_2 \leq C \Delta^{p-1} \int_{t_{k-1}}^{t_k}
\biggl(
\mathbb{E} [\vert X_s-X_{t_{k-1}} \vert^p |{\mathcal{F}}_{t_{k-1}}]+\vert X_{t_{k-1}} \vert^p
+ \mathbb{E} [\vert X_s\vert^p ]\biggr) ds.
\end{align*}
Therefore, using (\ref{estimateA1}), we get
$$
A_2 \leq C_T \Delta^p (1+\vert X_{t_{k-1}} \vert)^p,
$$
which concludes the proof of the Lemma.
\end{proof}

\vskip 12pt
The following lemma is the conditional version of~\cite[Theorem 11]{Guerin:02}.
\begin{lem} \label{sup}
For any $p>1$, $m \geq 1$ and $k\in \{1,...,N\}$, there exists a finite constant
$C_T$ such that, for $1 \leq i,l_1,...,l_m \leq d$,
\begin{align} \label{supremum} \nonumber
& \sup_{r_1,...,r_m,s \in [t_{k-1},t_k]} \mathbb{E} \biggl[ \int_0^1
\cdots\int_0^1|D^{l_1}_{(r_1,z_1)}\cdots D^{l_m}_{(r_m,z_m)}(X^i_s) |^p dz_1 \cdots
dz_m| {\mathcal{F}}_{t_{k-1}} \biggr] \\
& \qquad \qquad \qquad \qquad \leq C_T (1+\vert X_{t_{k-1}}
\vert)^p .
\end{align}
\end{lem}

\begin{proof}
We proceed by induction on $m$.
Suppose $m=1$. Let $z \in [0,1]$. For $r,s  \in
[t_{k-1},t_k]$ and $1 \leq i,l \leq d$,
we consider the stochastic
differential equation satisfied by the derivative
(cf.~\cite[Theorem 11]{Guerin:02})
\begin{align} \label{sde} \nonumber
D^{l}_{(r,z)}(X_s^i)=\sigma_{il}(X_r-Y_r(z))&+\int_r^s \int_0^1
\sum_{j,n=1}^d \partial_n \sigma_{ij}(X_u-Y_u(\alpha))
D^l_{(r,z)}(X_u^n) W^j(d\alpha,du) \\
&+\int_r^s \int_0^1 \sum_{n=1}^d \partial_n b_i(X_u-Y_u(\alpha))
D^l_{(r,z)}(X_u^n) d\alpha du.
\end{align}
Note that
$$
\sum_{i=1}^d \mathbb{E} \biggl[ \int_0^1
|D^{l}_{(r,z)}(X^i_s) |^p dz| {\mathcal{F}}_{t_{k-1}} \biggr]
\leq \sum_{i=1}^d  3^{p-1} (A_1+A_2+A_3),
$$
where
\begin{align*}
A_1&=\mathbb{E} \biggl[ \int_0^1
|\sigma_{il}(X_r-Y_r(z)) |^p dz| {\mathcal{F}}_{t_{k-1}} \biggr] \\
A_2&=\mathbb{E} \biggl[ \int_0^1
\biggl( \int_r^s \int_0^1
\sum_{j,n=1}^d \partial_n \sigma_{ij}(X_u-Y_u(\alpha))
D^l_{(r,z)}(X_u^n) W^j(d\alpha,du) \biggr)^p dz| {\mathcal{F}}_{t_{k-1}} \biggr] \\
A_3&=\mathbb{E} \biggl[ \int_0^1
\biggl( \int_r^s \int_0^1 \sum_{n=1}^d \partial_n b_i(X_u-Y_u(\alpha))
D^l_{(r,z)}(X_u^n) d\alpha du \biggr)^p dz | {\mathcal{F}}_{t_{k-1}} \biggr].
\end{align*}
Now, from the Lipschitz property of $\sigma$ and (\ref{estimateA1}), we have that
\begin{align*}
A_1 &\leq C_T (\mathbb{E}[\vert X_r- X_{t_{k-1}} \vert^p |
{\mathcal{F}}_{t_{k-1}}] + 1+ \vert X_{t_{k-1}}
\vert^p) \\
&\leq C_T (1+\vert X_{t_{k-1}}
\vert)^p.
\end{align*}
Moreover, using the bounds of the derivatives of $\sigma$,
Burkholder's and H\"older's inequalities for conditional expectations, it yields
$$
A_2 \leq C_T \mathbb{E} \biggl[ \int_0^1 \int_r^s \sum_{n=1}^d \vert D^l_{(r,z)}(X_u^n) \vert^p du dz| {\mathcal{F}}_{t_{k-1}} \biggr].
$$
Finally, the bounds of the derivatives of $b$ and H\"older's inequality imply that
$$
A_3 \leq  C_T \mathbb{E} \biggl[ \int_0^1 \int_r^s \sum_{n=1}^d \vert D^l_{(r,z)}(X_u^n) \vert^p du dz| {\mathcal{F}}_{t_{k-1}} \biggr].
$$
Hence, using Gronwall's Lemma, we conclude that
$$
\sum_{i=1}^d \mathbb{E} \biggl[ \int_0^1
|D^{l}_{(r,z)}(X^i_s) |^p dz| {\mathcal{F}}_{t_{k-1}} \biggr]
\leq C_T (1+\vert X_{t_{k-1}}
\vert)^p,
$$
which proves (\ref{supremum}) for $m=1$.

For $m>1$, consider the
stochastic differential equation satisfied by the iterated
derivative, for $r_1,...,r_m,s \in [t_{k-1}, t_k]$, $z_1,...,z_m \in [0,1]$, $1 \leq i, l_1,...l_m \leq d$,
\begin{align} \label{iterated}
 \nonumber &D^{l_1}_{(r_1,z_1)} \cdots D^{l_m}_{(r_m,z_m)} (X_s^i)\\  \nonumber
&=\sum_{n=1}^m D^{l_1}_{(r_1,z_1)} \cdots
D^{l_{n-1}}_{(r_{n-1},z_{n-1})} D^{l_{n+1}}_{(r_{n+1},z_{n+1})}
\cdots D^{l_m}_{(r_m,z_m)}(\sigma_{i l_n}(X_{r_n}-Y_{r_n}(z_n)) ) \\ \nonumber
&\qquad   + \sum_{j=1}^d \int_{r_1}^s \int_0^1 \cdots
\int_{r_m}^s \int_0^1
D^{l_1}_{(r_1,z_1)} \cdots D^{l_m}_{(r_m,z_m)} (\sigma_{ij}(X_u-Y_u(\alpha))) \, W^j(d\alpha,du) \\
& \qquad  +\int_{r_1}^s \int_0^1 \cdots \int_{r_m}^s \int_0^1
D^{l_1}_{(r_1,z_1)} \cdots D^{l_m}_{(r_m,z_m)}
(b_i(X_u-Y_u(\alpha))) \, d\alpha du.
\end{align}
Then, using the induction hypothesis and Gronwall's Lemma, one completes the desired proof.
\end{proof}

\vskip 12pt
The next result gives an upper bound for the derivative of $J_k+\Gamma_k$.
\begin{lem} \label{derivativegamma}
For any $p>1$ and $m \geq 1$, there exists a finite constant
$C_T>0$ such that, for all $i\in \{1,...,d\}$ and $k\in \{1,...,N\}$,
\begin{equation*}
(\mathbb{E} [ \Vert D^{(m)} (J^i_k+\Gamma_k^i)
\Vert^p_{\mathcal{H}_{t_{k-1}}^{\otimes m}}|
{\mathcal{F}}_{t_{k-1}}])^{1/p} \leq C_T \Delta^{1/2} (1+\vert X_{t_{k-1}}
\vert).
\end{equation*}
\end{lem}

\begin{proof}
Let $(r,z) \in [0,t] \times [0,1]$. Note that, for $i,l=1,...,d$,
\be \label{dersigma}
D^{l}_{(r,z)}(J^i_k)=\sigma_{i,l}(X_{t_{k-1}}-Y_{t_{k-1}}(z)) \,
1_{[t_{k-1},t_k]}(r), \ee
 and, therefore, the iterated derivative
$D^{(m)}_{(r,z)}(J^i_k)$ equals zero for $m>1$.

Hence, using the Lipschitz continuity of $\sigma$, we get
\begin{align*}
\mathbb{E} [ \Vert D^{(m)} (J^i_k)
\Vert^p_{\mathcal{H}_{t_{k-1}}^{\otimes m}}|
{\mathcal{F}}_{t_{k-1}}]&= \mathbb{E}
\biggl[\biggl(\int_{t_{k-1}}^{t_k}\int_0^1
\sum_{j=1}^d \vert \sigma_{ij}(X_{t_{k-1}}-Y_{t_{k-1}}(z))\vert^2 dr dz \biggr)^{p/2} | {\mathcal{F}}_{t_{k-1}} \biggr]\\
&\leq C_T \Delta^{p/2} (1+\vert X_{t_{k-1}} \vert)^p.
\end{align*}

On the other hand, for $r \in
[t_{k-1},t_k]$, and $1 \leq i,l \leq d$,
\begin{align*}
& \qquad D^{l}_{r,z}(\Gamma^i_k)=\sigma_{il}(X_r-Y_r(z))-
\sigma_{il}(X_{t_{k-1}}-Y_{t_{k-1}}(z)) \\
&+ \int_{r}^{t_k} \int_0^1 \sum_{j=1}^d
D_{(r,z)}^{l}(\sigma_{ij}(X_s-Y_s(\alpha))) \, W^j(d\alpha, ds)
+\int_{r}^{t_k} \int_0^1 D_{(r,z)}^{l}(b_i(X_s-Y_s(\alpha))) \,
d\alpha ds,
\end{align*}
and is equal to zero elsewhere. Therefore,
\begin{equation} \label{3}
\mathbb{E} [\Vert
D(\Gamma_k^i)\Vert_{\mathcal{H}_{t_{k-1}}}^p|{\mathcal{F}}_{t_{k-1}}]
\leq 3^{p-1}(A_1+A_2+A_3),
\end{equation}
where
\begin{align*}
A_1 &:=\mathbb{E}\biggl[\biggl( \int_{t_{k-1}}^{t_k}\int_0^1
\sum_{j=1}^d \vert \sigma_{ij}(X_r-Y_r(z))-
\sigma_{ij}(X_{t_{k-1}}-Y_{t_{k-1}}(z))\vert^2 dr dz \biggr)^{p/2}
|
{\mathcal{F}}_{t_{k-1}}\biggr], \\
A_2 &:= \mathbb{E}\biggl[\biggl( \int_{t_{k-1}}^{t_k}\int_0^1
\sum_{l=1}^d ( \int_{r}^{t_k} \int_0^1 \sum_{j=1}^d
D_{(r,z)}^{l}(\sigma_{ij}(X_s-Y_s(\alpha))) \, W^j(d\alpha,
ds))^2 dr
dz \biggr)^{p/2} | {\mathcal{F}}_{t_{k-1}}\biggr] \\
A_3 &:= \mathbb{E}\biggl[\biggl( \int_{t_{k-1}}^{t_k}\int_0^1
\sum_{l=1}^d (\int_{r}^{t_k} \int_0^1
D_{(r,z)}^{l}(b_i(X_s-Y_s(\alpha))) d\alpha ds)^2 dr dz
\biggr)^{p/2} | {\mathcal{F}}_{t_{k-1}}\biggr].
\end{align*}

From the proof of Lemma~\ref{p1derivativegamma} we get
$$
A_1 \leq C_T \Delta^p (1+\vert X_{t_{k-1}} \vert)^p.
$$
For the second term, use Burkholder's and H\"older's
inequalities for conditional expectations, the bounds of the
derivatives of $\sigma$ and Lemma~\ref{sup} to conclude that
\begin{align*}
A_2 &\leq C_T \Delta^p \sum_{l=1}^d \sup_{r,s \in [t_{k-1},t_k]}
\mathbb{E} \biggl[ \int_0^1 |D_{(r,z)}^{l}(X^i_s) |^p dz|
{\mathcal{F}}_{t_{k-1}} \biggr] \\
&\leq C_T \Delta^p (1+\vert X_{t_{k-1}} \vert)^p.
\end{align*}
Finally, using H\"older's inequality, the bounds for the derivative of
$b$ and Lemma~\ref{sup}, we obtain
\begin{eqnarray*}
A_3 &\leq& C_T \Delta^p \sum_{l=1}^d \sup_{r,s \in [t_{k-1},t_k]}
\mathbb{E}\biggl[ \int_0^1 |D_{(r,z)}^{l}(X^i_s) |^p dz|
{\mathcal{F}}_{t_{k-1}}\biggr] \\
&\leq& C_T \Delta^p (1+\vert X_{t_{k-1}} \vert)^p.
\end{eqnarray*}
Using (\ref{3}), it yields
\begin{equation} \label{new}
\mathbb{E} [\Vert
D(\Gamma_k^i)\Vert_{\mathcal{H}_{t_{k-1}}}^p|{\mathcal{F}}_{t_{k-1}}] \leq
C_T \Delta^p (1+\vert X_{t_{k-1}} \vert)^p.
\end{equation}

In order to treat the other derivatives we use the stochastic
differential equation satisfied by the iterated derivatives and
similar arguments to conclude that, for $m \geq 1$,
\begin{equation} \label{deriga}
\mathbb{E} [ \Vert D^{(m)}(\Gamma_k^i)
\Vert^p_{\mathcal{H}_{t_{m-1}}^{\otimes m}}|
{\mathcal{F}}_{t_{k-1}}] \leq C_T \Delta^p (1+\vert X_{t_{k-1}}
\vert)^p,
\end{equation}
which proves the Lemma.
\end{proof}

\vskip 12pt

As a consequence of Lemma~\ref{p1derivativegamma} and (\ref{deriga}) we obtain the following
estimate for the Sobolev norm of $\Gamma_k$.
\begin{cor} \label{SobolevGamma}
For any $p>1$ and $m\geq 0$, there exists a finite constant $C_T$ such that,
for $i\in \{1,...,d\}$ and $k\in \{1,...,N\}$,
$$
\Vert\Gamma_k^i\Vert_{m,p,t_{k-1}} \leq
C_T \Delta (1+\vert X_{t_{k-1}} \vert).
$$
\end{cor}

\vskip 12pt

We will also need the following lower bound for the determinant of the
Malliavin matrix of $J_k+\Gamma_k$.
\begin{lem} \label{determinant}
For any $p>1$ and $q>d$, there exists a finite constant
$C_T>0$ such that, for any $i\in \{1,...,d\}$, $k\in \{1,...,N\}$ and $0< \rho \leq 1$,
$$
\mathbb{E} [ (\textnormal{det} \, \gamma_{J_k+\rho \Gamma_k}
(t_{k-1}) )^{-p} | {\mathcal{F}}_{t_{k-1}} ]^{1/p} \leq  C_T
\Delta^{-d} (1+\vert X_{t_{k-1}} \vert)^{2q}.
$$
\end{lem}

\begin{proof}
In order to simplify the notation we write
$\gamma_k:=\gamma_{J_k+\rho \Gamma_k} (t_{k-1})$. Note that
$$
(\text{det} \, \gamma_k)^{1/d} \geq \inf_{\xi \in
\mathbb{R}^d,
\vert \xi \vert=1} \langle \gamma_k \xi, \xi \rangle,
$$
where
$$
\langle \gamma_k \xi, \xi \rangle =\sum_{l=1}^d
\int^{t_k}_{t_{k-1}} \int_0^1 \vert \sum_{i=1}^d D^{l}_{(r,z)}
(J_k^i+\rho \Gamma_k^i) \xi_i \vert^2 dz dr.
$$
Now, fix $h \in (0,1]$. Using the inequality $(a+b)^2\geq \frac{1}{2} a^2 -b^2$, we obtain that
\begin{align*}
\langle \gamma_k \xi, \xi \rangle &\geq \sum_{l=1}^d
\int^{t_k}_{t_k-h(t_k-t_{k-1})} \int_0^1 \vert \sum_{i=1}^d D^{l}_{(r,z)}
(J_k^i+\rho \Gamma_k^i) \xi_i \vert^2 dz dr \\
&\geq \sum_{l=1}^d
\int^{t_k}_{t_k-h(t_k-t_{k-1})} \int_0^1 \biggl( \frac{1}{2} (\sum_{i=1}^d D^{l}_{(r,z)}(J_k^i) \xi_i)^2
-(\sum_{i=1}^d D^{l}_{(r,z)} (\rho \Gamma_k^i) \xi_i)^2 \biggr) dz dr.
\end{align*}
Moreover, by (\ref{dersigma}) and (\ref{ev.min}), it yields
$$
\inf_{\xi \in \mathbb{R}^d, \vert \xi \vert=1} \langle \gamma_k
\xi, \xi \rangle \geq \frac{\lambda_1}{2} h \Delta - \sup_{\xi
\in \mathbb{R}^d, \vert \xi \vert=1} I_h,
$$
where
$$
I_h:=\sum_{l=1}^d \int^{t_k}_{t_k-h \Delta} \int_0^1
\biggl(\sum_{i=1}^d D^{l}_{(r,z)} (\rho \Gamma_k^i) \xi_i
\biggr)^2 dz dr.
$$

Using (\ref{new}), for $q>1$, we have that
$$
\mathbb{E} \biggl[ \sup_{\xi \in \mathbb{R}^d, \vert \xi \vert=1}
\vert I_h \vert ^q | {\mathcal{F}}_{t_{k-1}} \biggr] \leq C_T
h^{2q} \Delta^{2q} (1+\vert X_{t_{k-1}} \vert)^{2q}.
$$

We now use Proposition~\ref{deter} with $\epsilon_0=\Delta$,
$\alpha_1=1$, $\alpha_2=2$, $c_1=\frac{\lambda_1}{2}$, $c_2=C_T$,
$Z=\inf_{\xi \in \mathbb{R}^d, \vert \xi \vert=1} \langle \gamma_k
\xi, \xi \rangle$, $\epsilon=h \Delta$,
$X(\epsilon)=\frac{\lambda_1}{2} h \Delta$, $Y(\epsilon)=\sup_{\xi
\in \mathbb{R}^d, \vert \xi \vert=1} I_h$, $s=t_{k-1}$ and
$G_{t_{k-1}}=(1+\vert X_{t_{k-1}} \vert)^2$. Then, we obtain that
for any $q>d$,
\begin{align*}
\mathbb{E} [ (\text{det} \, \gamma_k )^{-p}|
{\mathcal{F}}_{t_{k-1}}]^{1/p}&\leq E[(\inf_{\xi \in \mathbb{R}^d,
\vert \xi \vert=1} \langle \gamma_k \xi, \xi \rangle)^{-dp}|
{\mathcal{F}}_{t_{k-1}}]^{1/p}
\\\nonumber
&\leq C_T \Delta^{-d} (1+\vert X_{t_{k-1}} \vert)^{2q},
\end{align*}
which concludes the desired result.
\end{proof}

\vskip 12pt
The next result gives an upper bound for the second term in
(\ref{meanvalue}).
\begin{prop} \label{upperboundderi}
There exists a constant $C_2>0$ depending only on $T$ and
independent of $k$ such that, for any $0< \rho \leq 1$, $z \in \mathbb{R}^d$ and $k\in \{1,...,N\}$, a.s.,
$$
\mathbb{E} \biggl[ \partial_{x^i} \phi_{\eta}
(X_{t_{k-1}}+J_k-z+\rho \Gamma_k) \Gamma^i_k |
{\mathcal{F}}_{t_{k-1}} \biggr]  \leq C_2 \Delta^{1/2-d/2}(1+\vert
X_{t_{k-1}} \vert)^D,
$$
where $D$ is polynomial of degree $3$ on $d$.
\end{prop}

\begin{proof}
Define
$$
\Phi_{\eta}(x)=\int_{-\infty}^{x_1} \cdots  \int_{-\infty}^{x_d}
\phi_{\eta}(u) du, \; \; x \in \mathbb{R}^d,
$$
and remark that
$$
\partial_{x^i} \phi_{\eta} (X_{t_{k-1}}+J_k-z+\rho
\Gamma_k)={\frac{\partial^{d+1} \Phi_{\eta}}{\partial x^i \partial
x^1 \cdots \partial x^d }}(X_{t_{k-1}}+J_k-z+\rho \Gamma_k).
$$
Using the version of the integration by parts formula given in Proposition~\ref{ipfc},
\begin{eqnarray*}
&&\mathbb{E} \biggl[ \partial_{x^i} \phi_{\eta}
(X_{t_{k-1}}+J_k-z+\rho \Gamma_k) \Gamma^i_k |
{\mathcal{F}}_{t_{k-1}} \biggr] \\
&&\qquad \qquad= \mathbb{E}\biggl[\Phi_{\eta}(X_{t_{k-1}}+J_k-z+\rho
\Gamma_k) H_{(1,...,d,i)}(J_k+\rho \Gamma_k, \Gamma^i_k ) |
{\mathcal{F}}_{t_{k-1}} \biggr].
\end{eqnarray*}
As $\int \phi_{\eta} =1$, by the Cauchy-Schwarz inequality, we
obtain
\begin{equation*}
 \mathbb{E} \biggl[\partial_{x^i} \phi_{\eta}(X_{t_{k-1}}+J_k-z+\rho \Gamma_k)
\Gamma^i_k | {\mathcal{F}}_{t_{k-1}} \biggr]  \leq \Vert
H_{(1,...,d,i)}(J_k+\rho \Gamma_k, \Gamma^i_k)
\Vert_{0,2,t_{k-1}}.
\end{equation*}
We now apply Proposition~\ref{normHc} with $\alpha=(1,...,d,i)$, $F=J_k+\rho\Gamma_k$
and $G=\Gamma^i_k$. For this, we use Lemma~\ref{derivativegamma} to prove (\ref{a1})
of Proposition~\ref{normHc} with $Z_{t_{k-1}}=\Delta^{1/2} (1+\vert X_{t_{k-1}} \vert)$,
and Lemma~\ref{determinant} with $q=d +\frac{1}{2}$ to prove (\ref{a2}) with $Y_{t_{k-1}}=(1+\vert X_{t_{k-1}} \vert)^{4d+1}$. Then, using Corollary~\ref{SobolevGamma}, we conclude that
\begin{align*}
&\mathbb{E} \biggl[\partial_{x^i}
\phi_{\eta}(X_{t_{k-1}}+J_k-z+\rho \Gamma_k) \Gamma^i_k |
{\mathcal{F}}_{t_{k-1}} \biggr]  \\&\leq C_T \Vert \Gamma_k^i
\Vert_{d+1, 2^{d+2}, t_{k-1}} \Delta^{-(d+1)/2} (1+\vert X_{t_{k-1}} \vert
)^{-(d+1)} \prod_{i=1}^{d+1} \sum_{j=1}^{i+1} (1+\vert X_{t_{k-1}} \vert)^{j(4d+1)} \\
&\leq C_T \Delta^{1/2-d/2} (1+\vert X_{t_{k-1}} \vert)^D,
\end{align*}
where $D$ is polynomial of degree $3$ in $d$. This proves the
desired bound.
\end{proof}

\vskip 12pt

Applying the bounds obtained in
Propositions~\ref{previouslowerbound} and~\ref{upperboundderi}
into (\ref{meanvalue}) we obtain the following lower bound for the
conditional density of $X_{t_k}$ given ${\mathcal{F}}_{t_{k-1}}$.
\begin{cor} \label{l}
Assume $0<\eta\leq \sqrt{\lambda_1 \Delta} $, and fix $z \in
\mathbb{R}^d$. Then, for almost all $(w,z)$ such that $ \vert
X_{t_{k-1}}(\omega)-z \vert \leq \sqrt{\lambda_1 \Delta}, $ it
holds
\begin{eqnarray*}
\mathbb{E} [\phi_{\eta}
(X_{t_{k}}-z)|{\mathcal{F}}_{t_{k-1}}] \geq \frac{1}{C_1
\Delta^{d/2} (1+\vert X_{t_{k-1}}\vert)^{d}}-C_2
\Delta^{1/2-d/2}(1+\vert X_{t_{k-1}} \vert)^D,
\end{eqnarray*}
where $C_1$, $C_2$ and $D$ are the constants obtained in
Propositions~\ref{previouslowerbound} and~\ref{upperboundderi}.
\end{cor}

\subsection{Proof of the lower bound}

We now fix $v \in \mathbb{R}^d$. Fix $x_0=X_0$, and let  $x_1$,
..., $x_{N-1}$, $x_N$ be $N$ $\mathcal{F}_0$-measurable points
defined by  $x_k = x_{k-1} +\frac{k-1}{N}( v-X_0)$ for $1\leq
k\leq N$. Remark that $x_N=v$, $|x_k|\leq |v-X_0|+|X_0|$, and
there exists a constant $C_3$ only depending on $\lambda_1$
and $T$, such that if $|x-x_k|\leq \frac{\sqrt{\lambda_1 T}}{2}$,
($x\in \mathbb{R}^d$), then
\begin{equation} \label{C3}
1+\vert x\vert\leq
C_3(1+|X_0|+|v-X_0|).
\end{equation}
 We choose the discretization size $N$ as
the smallest integer such that
\begin{equation*}
N \geq \frac{16 \vert v-X_0
\vert^2}{\lambda_1 t} +\frac{t}{M}+1,
\end{equation*}
where $$M=\frac{1}{ (2C_1C_2C_3^{d+D})^2(1+\vert X_0 \vert+\vert
X_0 - v \vert)^{2(d+D)}}.$$ The constants $C_1$, $C_2$ and $D$ are
defined in Propositions~\ref{previouslowerbound}
and~\ref{upperboundderi}.

This choice of $N$ will be justified by the computations below.
Note that, in particular, it implies that
\begin{equation*}
\frac{t}{N}=\Delta \leq M,
\end{equation*}
and that for each $1\leq k\leq N$,
\begin{equation}  \label{x}
\vert x_k-x_{k-1} \vert \leq \frac{\sqrt{\lambda_1 \Delta}}{4}.
\end{equation}

\vskip 12pt

We introduce the following sets, for $k=1,...,N$,
$$
A_k=  \{ \omega : \vert X_{t_{i-1}}(\omega)-x_i \vert  \leq
\frac{\sqrt{\lambda_1 \Delta}}{2}, i=1,...,k \} \in
\mathcal{F}_{t_{k-1}}.
$$

\begin{prop} \label{l1}
 Assume $0<\eta\leq \sqrt{\lambda_1 \Delta} $. Let $k\in
\{1,...,N\}$ and  consider $z \in \mathbb{R}^d$
 such that $ \vert x_k-z \vert \leq \frac{\sqrt{\lambda_1
\Delta}}{2}$. Then, a.s.
\begin{equation*}
\mathbb{E} [\phi_{\eta} (X_{t_{k}}-z)|{\mathcal{F}}_{t_{k-1}}]
 \geq \frac{1}{2C_1C_3^d \Delta^{d/2} (1+|X_0|+ |v-X_0|)^{d}}\,{\bf 1}_{A_k}.
\end{equation*}
\end{prop}

\begin{proof} Remark that if $\omega\in A_k$ and $ \vert x_k-z \vert \leq
\frac{\sqrt{\lambda_1 \Delta}}{2}$, then $ \vert
X_{t_{k-1}}(\omega)-z \vert \leq \sqrt{\lambda_1 \Delta}$.
Therefore, using Corollary~\ref{l}, (\ref{C3}), and the choice of $\Delta$,
we get
\begin{align*}
\mathbb{E} [\phi_{\eta} (X_{t_{k-1}}-z) | {\mathcal{F}}_{t_{k-1}}]
&\geq \frac{1}{C_1 \Delta^{d/2} (1+\vert X_{t_{k-1}}
\vert)^{d}}-C_2 \Delta^{1/2-d/2}(1+\vert
X_{t_{k-1}} \vert)^D \\
&\geq  \frac{1}{\Delta^{d/2}} \frac{1}{2C_1C_3^d (1+|X_0|+
|v-X_0|)^{d}}.
\end{align*}
\end{proof}

\begin{prop} \label{Pa}
There exists a constant $C_4>0$ only depending on $\lambda_1$,
$\lambda_2$ and $T$ such that, for any $k\in \{1,...,N\}$,
$$
\mathbb{P}_{X_0}(A_k) \geq \frac{1}{C_4 (1+\vert
X_0\vert+|v-X_0|)^{d}} \mathbb{P}_{X_0}(A_{k-1}).
$$
\end{prop}

\begin{proof}
Let $0 < \eta < \sqrt{\lambda_1 \Delta}$. As $A_k=A_{k-1} \cap \{
|X_{t_{k-1}}-x_k| \leq \frac{\sqrt{\lambda_1 \Delta}}{2}\}$ and
using the fact that $\int \phi_{\eta}=1$, we have
\begin{eqnarray*}
&&\mathbb{P}_{X_0}(A_k) \\
&=& \mathbb{E}_{X_0}[1_{A_{k-1}} \mathbb{E}[1_{\{
|X_{t_{k-1}}-x_k| \leq \frac{\sqrt{\lambda_1\Delta}}{2}\}} |
\mathcal{F}_{t_{k-2}}]] \\
&=&\mathbb{E}_{X_0}[1_{A_{k-1}} \int_{\mathbb{R}^d} \mathbb{E}[
\phi_{\eta} (X_{t_{k-1}}-z) 1_{\{\vert X_{t_{k-1}}-x_k \vert \leq
\frac{\sqrt{\lambda_1 \Delta}}{2} \}} |
\mathcal{F}_{t_{k-2}}]dz ] \\
&\geq &\mathbb{E}_{X_0} \biggl[1_{A_{k-1}} \int_{\vert z-x_{k-1}
\vert \leq \sqrt{\lambda_1 \Delta}/4-\eta} \mathbb{E}[ \phi_{\eta}
(X_{t_{k-1}}-z) 1_{\{\vert X_{t_{k-1}}-x_k \vert \leq
\frac{\sqrt{\lambda_1 \Delta}}{2} \}} |
\mathcal{F}_{t_{k-2}}]dz  \biggr] \\
&=& \mathbb{E}_{X_0} \biggl[1_{A_{k-1}} \int_{\vert z-x_{k-1}
\vert \leq \sqrt{\lambda_1 \Delta}/4-\eta} \mathbb{E}[ \phi_{\eta}
(X_{t_{k-1}}-z) | \mathcal{F}_{t_{k-2}}]dz \biggr].
\end{eqnarray*}
The last equality follows from (\ref{x}) and the fact that
\begin{eqnarray*}
\vert X_{t_{k-1}}-x_k\vert &\leq&  \vert X_{t_{k-1}}-z\vert+\vert
z-x_{k-1} \vert +\vert
x_{k-1}-x_k \vert \\
&\leq& \eta + \frac{\sqrt{\lambda_1 \Delta}}{4}-\eta
+\frac{\sqrt{\lambda_1 \Delta}}{4}=\frac{\sqrt{\lambda_1
\Delta}}{2}.
\end{eqnarray*}
Take $\eta=\frac{\sqrt{\lambda_1 \Delta}}{8}$. Using
Proposition~\ref{l1} we obtain
\begin{eqnarray*}
\mathbb{P}_{X_0}(A_k) &\geq& \mathbb{E}_{X_0} \biggl[1_{A_{k-1}}
\int_{\vert z-x_{k-1}\vert \leq \sqrt{\lambda_1 \Delta}/8}
\mathbb{E}[ \phi_{\eta} (X_{t_{k-1}}-z)
| \mathcal{F}_{t_{k-2}}]dz \biggr] \\
&\geq& \mathbb{E}_{X_0} \biggl[1_{A_{k-1}} \int_{\vert z-x_{k-1}
\vert \leq \sqrt{\lambda_1 \Delta}/8}\frac{1}{2C_1C_3^d
\Delta^{d/2} (1+|X_0|+ |v-X_0|)^{d}}\,  dz \biggr] \\
&\geq& \frac{1}{2 C_1C_3^d \Delta^{d/2} (1+\vert X_0
\vert+|v-X_0|)^{d}} \biggl( \frac{\sqrt{\lambda_1 \Delta}}{8}
\biggr)^d \mathbb{P}_{X_0}(A_{k-1}).
\end{eqnarray*}
This concludes the proof of the Proposition.
\end{proof}
\vskip 12pt

We now conclude the proof of the lower bound. Let us apply
Proposition~\ref{l1} with $k=N$ and $z=v$ and an iteration of
Proposition \ref{Pa}.
\begin{align*}
\mathbb{E}[\phi_{\eta}(X_{t_N}-v) | X_0] &\geq \mathbb{E}[
\mathbb{E} [\phi_{\eta}(X_{t_N}-v) | \mathcal{F}_{t_{N-1}}]
1_{A_N} | X_0] \\
&\geq \frac{1}{2 C_1C_3^d \Delta^{d/2} (1+\vert X_0
\vert+|v-X_0|)^{d}} \mathbb{P}_{X_0}(A_N) \\
&\geq \frac{C_4}{2 C_1C_3^d  }\frac{N^{d/2}}{t^{d/2} }
\biggl(\frac{1}{C_4(1+\vert X_0
\vert+\vert v- X_0 \vert)^{d} }\biggr)^N \mathbb{P}_{X_0}(A_1).
\end{align*}
The choice of $N$ implies that $\mathbb{P}_{X_0}(A_1)=1$ a.s., and that
$$
\frac{16}{\lambda_1 t} \vert
v - X_0 \vert^2 +\frac{t}{M}+1 \leq N \leq
\frac{16}{\lambda_1 t} \vert v - X_0 \vert^2+\frac{t}{M}  + 2.
$$
Therefore, we obtain that
$$
\mathbb{E}[\phi_{\eta}(X_{t_N}-v) | X_0] \geq   \frac{1}{t^{d/2}
c_1(T,v,X_0)} e^{-c_2(T,v,X_0) \frac{\vert v-X_0 \vert^2}{t}},
$$
where the constants $c_1(T,v,X_0)$ and $c_2(T,v,X_0)$ can be explicitely
given as functions of $T, v, X_0, \lambda_1$ and $\lambda_2$.

\vskip 12pt

This concludes the proof of Theorem~\ref{maintheorem} (a).

\section{The upper bound}

In this section we prove Theorem~\ref{maintheorem} (b).

Let $T>0$, $0<t \leq 0$ and $v \in \mathbb{R}^d$ be fixed. Apply
Cauchy-Schwarz  inequality for conditional expectations to the
expression of Corollary~\ref{expression} with $\sigma=\{i \in \{
1,...,d\}: v_i \geq 0 \}$ to find that
\begin{equation} \label{upperbound}
f_{X_0}(t,v) \leq  (\mathbb{P}_{X_0} \{\vert X_t \vert \geq \vert
v \vert\})^{1/2} (\mathbb{E}_{X_0} [ (H^0_{(1,...,d)}(X_t,1))^2
])^{1/2},\quad P_0 \hbox{-a.s.}
\end{equation}
We estimate the first factor $\mathbb{P}_{X_0} \{\vert X_t \vert
\geq \vert v \vert\}^{1/2}$ using an exponential martingale
inequality. In order to deal with bounded coefficients, we
consider the SDE satisfied by a logarithmic transformation of our
process $X_t$. On the other hand, to obtain an upper bound for the
second factor $(\mathbb{E}_{X_0} [ (H^0_{(1,...,d)}(X_t,1))^2
])^{1/2}$ of order $t^{-d/2}$, we will use
Proposition~\ref{normHc} and precise estimates on the Sobolev
norms of $X_t$.

\vskip 12pt
This is given in the following two lemmas.
\begin{lem} \label{largedeviation}
There exist finite constants $c_1$ and $c_2$ only depending on $T$
such that for any $t \in (0, T]$ and $v \in \mathbb{R}^d$,
$P_0$-a.s.
\begin{equation*}
(\mathbb{P}_{X_0} \{\vert X_t \vert \geq \vert v \vert
\})^{1/2}\leq \exp \biggl(-\frac{( \ln(1+|v|^2)-\ln(1+|X_0|^2)-c_1
t )^2}{c_2 t} \biggr).
\end{equation*}
\end{lem}

\begin{proof}
Consider $Z_t=\ln(1+\vert X_t \vert^2)$. From the $d$-dimensional
It\^o's formula,
\begin{eqnarray*}
&&Z_t=\ln(1+\vert X_0 \vert^2)+\int_0^t \int_0^1 \sum_{i,j=1}^d
\frac{2 X_s^i}{1+\vert X_s \vert^2} \sigma_{ij}(X_s-Y_s(\alpha))
W^j(d \alpha, ds) \\
&& \qquad \qquad +\int_0^t \int_0^1 \sum_{i=1}^d \frac{2
X_s^i}{1+\vert X_s \vert^2} b_i(X_s-Y_s(\alpha)) d \alpha ds \\
&& \qquad \qquad +\int_0^t \int_0^1 \sum_{i,j=1}^d
\frac{1}{1+\vert X_s \vert^2} (\sigma_{ij}(X_s-Y_s(\alpha)))^2 d
\alpha ds \\
&& \qquad \qquad -\int_0^t \int_0^1 \sum_{i,j,k=1}^d \frac{2X_s^i
X_s^k}{(1+\vert X_s \vert^2)^2} \sigma_{ij}(X_s-Y_s(\alpha))
\sigma_{kj}(X_s-Y_s(\alpha))d \alpha ds.
\end{eqnarray*}
Using the Lipschitz property of $b$, we have that
\begin{align*}
\biggl| \int_0^t \int_0^1 \sum_{i=1}^d \frac{2 X_s^i \,
b_i(X_s-Y_s(\alpha))}{1+\vert X_s \vert^2} d \alpha ds  \biggr|
&\leq C( t +  t \mathbb{E}[\sup_{0 \leq s \leq T} \vert X_s \vert]) \\
& \leq C_1 t.
\end{align*}
Equally, from the Lipschitz property of $\sigma$,
$$
\biggl|\int_0^t \int_0^1 \sum_{i,j=1}^d \frac{1}{1+\vert X_s
\vert^2} (\sigma_{ij}(X_s-Y_s(\alpha)))^2 d \alpha ds \biggr| \leq
C_2 t,
$$
and
$$
\biggl|\int_0^t \int_0^1\sum_{i,j,k=1}^d \frac{2X_s^i
X_s^k}{(1+\vert X_s \vert^2)^2} \sigma_{ij}(X_s-Y_s(\alpha))
\sigma_{ij}(X_s-Y_s(\alpha)) d \alpha ds \biggr| \leq C_3 t.
$$
Hence, we obtain
\begin{align} \label{exponential} \nonumber
\mathbb{P}_{X_0} \{ \vert X_t \vert \geq \vert v \vert \} &\leq
\mathbb{P}_{X_0} \{Z_t \geq \ln (1+\vert v \vert^2)
 \} \\
&\leq  \mathbb{P}_{X_0} \{ M_t \geq \ln (1+\vert v
\vert^2)-\ln(1+\vert X_0 \vert^2)-c_1 t \},
\end{align}
where $c_1:=C_1+C_2+C_3$ and
$$M_t= \int_0^t \int_0^1 \sum_{i,j=1}^d
\frac{2 X_s^i}{1+\vert X_s \vert^2} \sigma_{ij}(X_s-Y_s(\alpha))
W^j(d \alpha, ds) $$ is a continuous martingale with respect to
$\mathcal{F}_t$ and with increasing process given by
$$\langle M \rangle_t=\int_0^t \int_0^1 \sum_{j=1}^d \biggl( \sum_{i=1}^d
\frac{2 X_s^i}{1+\vert X_s \vert^2} \sigma_{ij}(X_s-Y_s(\alpha))
\biggr)^2 d \alpha ds.$$ Again, using the Lipschitz property of
$\sigma$, we get that
$$\langle M \rangle_t \leq c t. $$
Finally, applying the exponential martingale inequality to
(\ref{exponential}), we obtain that $P_0$-a.s.
$$
\mathbb{P}_{X_0} \{ \vert X_t \vert \geq \vert v \vert \} \leq
\exp \biggl( - \frac{(\ln (1+\vert v \vert^2)-\ln(1+\vert X_0
\vert^2)-c_1 t)^2}{2 c t }\biggr).
$$
\end{proof}

\begin{lem} \label{upperboundH}
There exists a finite constant $c_3(T,X_0)>0$ such that $P_0$-a.s.
\begin{equation*}
(\mathbb{E}_{X_0} [ (H^0_{(1,...,d)}(X_t,1))^2 ])^{1/2} \leq
c_3(T,X_0) t^{-d/2},
\end{equation*}
for all $t \in (0,T]$.
\end{lem}

\begin{proof}
In order to prove this result, it suffices to prove that for any
$p>1$ and $m\geq 1$ there exist finite constants
$c_1(m,p,T,X_0)>0$ and $c_2(p,T,X_0) \geq 0$ such that
\begin{itemize}
\item[\textnormal{(i)}] $\mathbb{E}_{X_0}[\Vert D^{(m)}(X_t^i)
\Vert^p_{\mathcal{H}_0^{\otimes m}} ]^{1/p} \leq c_1(m,p,T,X_0) \,
t^{1/2}, \; i=1,...,d$; \vskip 12pt \item[\textnormal{(ii)}]
$\mathbb{E}_{X_0}[(\textnormal{det}
 \,\gamma_{X_t}(0) )^{-p}]^{1/p} \leq c_2(p,T,X_0) \, t^{-d}$.
\end{itemize}
Then, Proposition~\ref{normHc} with $s=0$ and $G=1$ concludes the
desired estimate.

\vskip 12pt We start proving (i). We proceed by induction on $m$.
For $m=1$, consider the stochastic differential equation
(\ref{sde}). Then, using H\"older's inequality for conditional
expectations, and Lemma \ref{sup}, we obtain,
\begin{align*}
\mathbb{E}_{X_0}[\Vert D(X_t^i) \Vert^p_{\mathcal{H}_0}]
&=\mathbb{E}_{X_0}\biggl[\biggl(\int_0^t \int_0^1 \vert
D_{(r,z)}(X_t^i)\vert^2
dr dz \biggr)^{p/2} \biggr] \\
&\leq t^{p/2} \biggl( \sup_{0\leq r \leq T}
\mathbb{E}_{X_0}\biggl[ \int_0^1 |D_{(r,z)}(X_t^i)|^p
dz\biggr]\biggr)\\
&\leq t^{p/2} C_T(1+|X_0|)^p.
\end{align*}
Then, the case $m>1$ follows along the same lines using the
stochastic differential equation satisfies by the iterated
derivative (\ref{iterated}) together with Lemma \ref{sup}.

We now prove (ii). Fix $\epsilon \in (0, 1/2]$ so that $t/2 \leq
t(1-\epsilon) < t$. From a similar argument as in Lemma
\ref{determinant}, it follows that
\begin{align*}
(\text{det} \, \gamma_{X_t}(0))^{1/d} &\geq \inf_{\xi \in
\mathbb{R}^d,
\vert \xi \vert=1} \langle \gamma_{X_t}(0)\xi, \xi \rangle \\
&\geq \frac{1}{2} m \tilde{c} t \epsilon - \sup_{\xi \in
\mathbb{R}^d, \vert \xi \vert=1} I_{\epsilon},
\end{align*}
where $m$ is defined in (\ref{boundh}), $\tilde{c}$ denotes the infimum of the function
\begin{equation*}
F(\xi,t)=\mathbb{E}[\vert X_t \vert^2 \vert \xi \vert^2- \langle
X_t, \xi \rangle^2]
\end{equation*}
on the compact set $\{ r\in
[\frac{t}{2},t]\} \times \{ \xi \in \mathbb{R}^d:\vert \xi
\vert=1\}$, and
\begin{eqnarray*}
&&I_{\epsilon}:=\sum_{k=1}^d \int_{t(1-\epsilon)}^t \int_0^1 \biggl(
\sum_{i=1}^d \xi_i \int_r^t \int_0^1 \sum_{j,l=1}^d
\partial_l \sigma_{ij}(X_s-Y_s(\alpha)) D^{k}_{(r,z)}(X^l_s)
W^j(d\alpha, ds) \\
&&\qquad \qquad\qquad+\sum_{i=1}^d \xi_i \int_r^t \int_0^1 \sum_{l=1}^d
\partial_l b_i(X_s-Y_s(\alpha)) D^k_{(r,z)}(X^l_s) d\alpha ds
\biggr)^2 dz dr.
\end{eqnarray*}

By some straightforward computations, using Burkholder's and
H\"older's inequalities and Lemma \ref{sup}, we obtain for any
$q>1$
\begin{eqnarray*}
\mathbb{E}_{X_0} \biggl[ \sup_{\xi \in \mathbb{R}^d, \vert \xi
\vert=1} |I_{\epsilon}|^q \biggr] &\leq& C_T (t \epsilon)^{2q}
\sup_{0\leq r,s \leq T}
\mathbb{E}_{X_0} \biggl[ \int_0^1 |D_{(r,z)} (X_s)|^{2q} dz  \biggr] \\
&\leq& C_T (1+|X_0|)^{2q}(t \epsilon)^{2q}.
\end{eqnarray*}
Consequently, applying Proposition~\ref{deter} with $Z=\inf_{\xi \in
\mathbb{R}^d,
\vert \xi \vert=1} \langle \gamma_{X_t}(0)\xi, \xi \rangle$, $\alpha_1=1$, $\alpha_2=2$ and $\epsilon_0=t$, we conclude that
$$
\mathbb{E}_{X_0}[(\textnormal{det}  \,\gamma_{X_t}(0)
)^{-p}]^{1/p} \leq  C(T,X_0) t^{-d},
$$
which proves (ii).

\end{proof}

\vskip 12pt

Substituting the results of Lemmas~\ref{largedeviation} and~\ref{upperboundH}
into the expression (\ref{upperbound}), we obtain that
$$
f_{X_0}(t,v) \leq c_3(T,X_0) \, t^{-d/2} e^{-\frac{(\ln(1+ \vert
v\vert^2)-\ln(1+ \vert X_0\vert^2)-c_1 t )^2}{c_2 t}}.
$$
 This concludes
the proof of the upper bound of Theorem~\ref{maintheorem}.

\Section{Acknowledgements}

The authors would like to thank V. Bally and A. Kohatsu-Higa for
all the fruitful discussions on the subject.

\end{document}